\documentclass[11pt,a4paper]{article}
\usepackage{iftex}
\ifPDFTeX
  \usepackage[utf8]{inputenc}
  \usepackage[T1]{fontenc}
\else
  \usepackage{fontspec}
  \IfFontExistsTF{cmunrm.otf}{\setmainfont{cmunrm.otf}[BoldFont=cmunbx.otf,ItalicFont=cmunti.otf,BoldItalicFont=cmunbi.otf]}{\setmainfont{DejaVu Serif}}
\fi
\usepackage[english]{babel}
\usepackage{amsmath,amssymb,amsthm,booktabs,array,adjustbox,needspace,longtable}
\usepackage[top=19mm,bottom=19mm,left=19mm,right=19mm,headsep=4mm]{geometry}
\usepackage[hidelinks]{hyperref}
\hypersetup{pdfauthor={Muratova Kh.A., Tyan M.D.},pdftitle={Derivations of Low-Dimensional Filiform Lie Superalgebras}}
\usepackage{fancyhdr}
\usepackage{titlesec}
\titlespacing*{\section}{0pt}{2ex plus .3ex}{.8ex}
\titlespacing*{\subsection}{0pt}{1.5ex plus .2ex}{.6ex}
\DeclareMathOperator{\Der}{Der}
\DeclareMathOperator{\rk}{rank}
\newtheorem{theorem}{Theorem}
\newtheorem{lemma}{Lemma}
\newtheorem{corollary}{Corollary}
\theoremstyle{definition}
\newtheorem{example}{Example}
\newcommand{\fitmath}[1]{\adjustbox{max width=\linewidth}{$\displaystyle #1$}}
\title{\bfseries Derivations of Low-Dimensional\\ Filiform Lie Superalgebras}
\author{\normalsize\bfseries Muratova Kh.A.,\quad Tyan M.D.\\[3pt]
{\small\mdseries\itshape Kimyo International University in Tashkent, Uzbekistan}\\[3pt]
\small E-mail: \href{mailto:muratovakhosiyat12@gmail.com}{muratovakhosiyat12@gmail.com},\\[-1pt]
\small \href{mailto:tyanmariya2006@gmail.com}{tyanmariya2006@gmail.com}}
\date{}
\begin{document}
\thispagestyle{plain}
\makeatletter
\begin{center}
{\LARGE\@title\par}
\vspace{.65em}
{\@author\par}
\end{center}
\makeatother
\vspace{-.4em}
\begin{center}
\begin{minipage}{0.94\linewidth}\small
\textbf{Abstract.} We determine the even and odd derivations of the 91 entries in Gilg's classification list of complex filiform Lie superalgebras of superdimensions $(2\mid2)$, $(2\mid3)$, $(2\mid4)$, $(3\mid2)$, $(3\mid3)$, $(4\mid2)$, $(5\mid2)$ and $(6\mid2)$. A unified linear system yields explicit matrix parametrizations and the dimensions of both homogeneous components of the derivation superalgebra. The results are given in block-matrix form.\\[3pt]
\end{minipage}
\end{center}

\section{Introduction}
Derivations of nilpotent Lie superalgebras play an important role in the study of their solvable extensions. In particular, the adjoint action of homogeneous elements of an extension on its nilradical defines homogeneous derivations. Methods for describing solvable Lie and Leibniz superalgebras with a prescribed nilradical have been developed in \cite{CamachoFernandezNavarro2020,CamachoNavarroOmirov2022}. Related questions for Leibniz superalgebras are considered in \cite{KhudoyberdiyevMuratova2024}.

A classification of low-dimensional filiform Lie superalgebras was given by Gilg \cite{Gilg2001}. Deformations of the model filiform Lie superalgebra were studied in \cite{Gilg2004,KhakimdjanovNavarro2010}, and further classification results were obtained in \cite{Navarro2016}. The cohomology of filiform Lie superalgebras was investigated in \cite{YangLiu2018}, while local derivations of solvable and model filiform superalgebras were studied in \cite{CamachoNavarroOmirov2023,ShengLiuLiu2024}.

The aim of this paper is to determine both homogeneous components of $\Der(F)=\Der_{\bar0}(F)\oplus\Der_{\bar1}(F)$ for eight families in \cite{Gilg2001}. We retain the original numbering of the cases. Rather than presenting repetitive case-by-case proofs, we employ a unified method that directly establishes the completeness of the matrix parametrizations. The explicit description of derivations provides a basis for further investigations of solvable extensions.

\section{Preliminaries}
Let $F=F_{\bar0}\oplus F_{\bar1}$ be a complex Lie superalgebra. For a homogeneous element $x$, denote its parity by $|x|\in\{0,1\}$. Set $C^0(F_{\bar\alpha})=F_{\bar\alpha}$ and $C^{k+1}(F_{\bar\alpha})=[F_{\bar0},C^k(F_{\bar\alpha})]$ for $\alpha=0,1$. In the notation of \cite{Gilg2001}, a filiform Lie superalgebra $F_{n,m}$ has superdimension $(n+1\mid m)$ and satisfies
\[
 C^{n}(F_{\bar0})=0,\quad C^{n-1}(F_{\bar0})\ne0,\qquad
 C^{m}(F_{\bar1})=0,\quad C^{m-1}(F_{\bar1})\ne0.
\]
Throughout this paper, we consider only nilpotent Lie superalgebras. We retain the notation $F_{n,m}^{(i)}$ for the $i$th case in \cite[pp.~469--477]{Gilg2001} and use the adapted basis
\begin{equation}\label{eq:basis}
\mathcal B=(e_1,\ldots,e_p,y_1,\ldots,y_q),\qquad p=n+1,\quad q=m,
\end{equation}
where $e_i=X_{i-1}$ and $y_j=Y_j$ in Gilg's notation. The following relations are common to all the multiplication tables considered:
\begin{equation}\label{eq:chains}
 [e_1,e_i]=e_{i+1}\quad (2\le i\le n),\qquad
 [e_1,y_j]=y_{j+1}\quad (1\le j<m).
\end{equation}
\subsection{Multiplication tables}
For convenience, we present the multiplication tables of all 91 cases under consideration from Gilg's classification \cite[pp.~469--477]{Gilg2001}, expressed in the basis \eqref{eq:basis}.
In each table, the brackets \eqref{eq:chains} common to the family are displayed above the column headings. Each row lists the additional nonzero brackets, grouped according to the parities of their arguments.
A dash indicates that there are no additional nonzero brackets of the corresponding type. All remaining products are obtained from the super-skew-symmetry relation
\[
 [x,y]=-(-1)^{|x||y|}[y,x]
\]
or are zero. Thus, the common relations together with the corresponding row completely determine the multiplication in $F_{n,m}^{(i)}$. Each case identifier also links to its matrices in Appendix~\ref{app:matrices}.
\begingroup
\small
\setlength{\tabcolsep}{3pt}
\renewcommand{\arraystretch}{1.04}
\setlength{\LTpre}{3pt}
\setlength{\LTpost}{4pt}
\setlength{\LTcapwidth}{\linewidth}

\begin{longtable}{@{}c>{\raggedright\arraybackslash}p{42mm}>{\raggedright\arraybackslash}p{42mm}>{\raggedright\arraybackslash}p{44mm}@{}}

\caption{Multiplication tables for the family $F_{1,2}$ (superdimension $(2\mid 2)$).}\label{tab:mult-1-2}\\
\toprule
\multicolumn{4}{@{}l@{}}{$[e_{1},y_{1}]=y_{2}$}\\[3pt]
\midrule
Case & $[F_{\bar0},F_{\bar0}]$ & $[F_{\bar0},F_{\bar1}]$ & $[F_{\bar1},F_{\bar1}]$\\
\midrule
\endfirsthead

\multicolumn{4}{@{}l@{}}{\textit{Table \thetable\ (continued; family $F_{1,2}$)}}\\[3pt]
\toprule
\multicolumn{4}{@{}l@{}}{$[e_{1},y_{1}]=y_{2}$}\\[3pt]
\midrule
Case & $[F_{\bar0},F_{\bar0}]$ & $[F_{\bar0},F_{\bar1}]$ & $[F_{\bar1},F_{\bar1}]$\\
\midrule
\endhead

\midrule
\multicolumn{4}{r@{}}{\textit{Continued on the next page}}\\
\endfoot
\bottomrule
\endlastfoot

\hyperref[case:1-2-1]{$F_{1,2}^{(1)}$} & \textemdash & \textemdash & \textemdash \\

\hyperref[case:1-2-2]{$F_{1,2}^{(2)}$} & \textemdash & $[e_{2},y_{1}]=y_{2}$ & \textemdash \\

\hyperref[case:1-2-3]{$F_{1,2}^{(3)}$} & \textemdash & \textemdash & $[y_{1},y_{1}]=e_{2}$ \\

\end{longtable}

\begin{longtable}{@{}c>{\raggedright\arraybackslash}p{42mm}>{\raggedright\arraybackslash}p{42mm}>{\raggedright\arraybackslash}p{44mm}@{}}

\caption{Multiplication tables for the family $F_{1,3}$ (superdimension $(2\mid 3)$).}\label{tab:mult-1-3}\\
\toprule
\multicolumn{4}{@{}l@{}}{$[e_{1},y_{1}]=y_{2},\quad [e_{1},y_{2}]=y_{3}$}\\[3pt]
\midrule
Case & $[F_{\bar0},F_{\bar0}]$ & $[F_{\bar0},F_{\bar1}]$ & $[F_{\bar1},F_{\bar1}]$\\
\midrule
\endfirsthead

\multicolumn{4}{@{}l@{}}{\textit{Table \thetable\ (continued; family $F_{1,3}$)}}\\[3pt]
\toprule
\multicolumn{4}{@{}l@{}}{$[e_{1},y_{1}]=y_{2},\quad [e_{1},y_{2}]=y_{3}$}\\[3pt]
\midrule
Case & $[F_{\bar0},F_{\bar0}]$ & $[F_{\bar0},F_{\bar1}]$ & $[F_{\bar1},F_{\bar1}]$\\
\midrule
\endhead

\midrule
\multicolumn{4}{r@{}}{\textit{Continued on the next page}}\\
\endfoot
\bottomrule
\endlastfoot

\hyperref[case:1-3-1]{$F_{1,3}^{(1)}$} & \textemdash & \textemdash & \textemdash \\

\hyperref[case:1-3-2]{$F_{1,3}^{(2)}$} & \textemdash & \textemdash & $[y_{1},y_{1}]=e_{2}$ \\

\hyperref[case:1-3-3]{$F_{1,3}^{(3)}$} & \textemdash & \textemdash & $[y_{2},y_{2}]=e_{2}$\newline $[y_{1},y_{3}]=-e_{2}$ \\

\hyperref[case:1-3-4]{$F_{1,3}^{(4)}$} & \textemdash & $[e_{2},y_{1}]=y_{2}$\newline $[e_{2},y_{2}]=y_{3}$ & \textemdash \\

\hyperref[case:1-3-5]{$F_{1,3}^{(5)}$} & \textemdash & $[e_{2},y_{1}]=y_{3}$ & \textemdash \\

\hyperref[case:1-3-6]{$F_{1,3}^{(6)}$} & \textemdash & $[e_{2},y_{1}]=y_{2}+y_{3}$\newline $[e_{2},y_{2}]=y_{3}$ & \textemdash \\

\end{longtable}

\begin{longtable}{@{}c>{\raggedright\arraybackslash}p{42mm}>{\raggedright\arraybackslash}p{42mm}>{\raggedright\arraybackslash}p{44mm}@{}}

\caption{Multiplication tables for the family $F_{1,4}$ (superdimension $(2\mid 4)$).}\label{tab:mult-1-4}\\
\toprule
\multicolumn{4}{@{}l@{}}{$[e_{1},y_{1}]=y_{2},\quad [e_{1},y_{2}]=y_{3},\quad [e_{1},y_{3}]=y_{4}$}\\[3pt]
\midrule
Case & $[F_{\bar0},F_{\bar0}]$ & $[F_{\bar0},F_{\bar1}]$ & $[F_{\bar1},F_{\bar1}]$\\
\midrule
\endfirsthead

\multicolumn{4}{@{}l@{}}{\textit{Table \thetable\ (continued; family $F_{1,4}$)}}\\[3pt]
\toprule
\multicolumn{4}{@{}l@{}}{$[e_{1},y_{1}]=y_{2},\quad [e_{1},y_{2}]=y_{3},\quad [e_{1},y_{3}]=y_{4}$}\\[3pt]
\midrule
Case & $[F_{\bar0},F_{\bar0}]$ & $[F_{\bar0},F_{\bar1}]$ & $[F_{\bar1},F_{\bar1}]$\\
\midrule
\endhead

\midrule
\multicolumn{4}{r@{}}{\textit{Continued on the next page}}\\
\endfoot
\bottomrule
\endlastfoot

\hyperref[case:1-4-1]{$F_{1,4}^{(1)}$} & \textemdash & \textemdash & \textemdash \\

\hyperref[case:1-4-2]{$F_{1,4}^{(2)}$} & \textemdash & \textemdash & $[y_{1},y_{1}]=e_{2}$ \\

\hyperref[case:1-4-3]{$F_{1,4}^{(3)}$} & \textemdash & \textemdash & $[y_{2},y_{2}]=e_{2}$\newline $[y_{1},y_{3}]=-e_{2}$ \\

\hyperref[case:1-4-4]{$F_{1,4}^{(4)}$} & \textemdash & $[e_{2},y_{1}]=y_{2}$\newline $[e_{2},y_{2}]=y_{3}$\newline $[e_{2},y_{3}]=y_{4}$ & \textemdash \\

\hyperref[case:1-4-5]{$F_{1,4}^{(5)}$} & \textemdash & $[e_{2},y_{1}]=y_{3}$\newline $[e_{2},y_{2}]=y_{4}$ & \textemdash \\

\hyperref[case:1-4-6]{$F_{1,4}^{(6)}$} & \textemdash & $[e_{2},y_{1}]=y_{4}$ & \textemdash \\

\hyperref[case:1-4-7]{$F_{1,4}^{(7)}$} & \textemdash & $[e_{2},y_{1}]=y_{2}+y_{3}$\newline $[e_{2},y_{2}]=y_{3}+y_{4}$\newline $[e_{2},y_{3}]=y_{4}$ & \textemdash \\

\hyperref[case:1-4-8]{$F_{1,4}^{(8)}$} & \textemdash & $[e_{2},y_{1}]=y_{2}+y_{4}$\newline $[e_{2},y_{2}]=y_{3}$\newline $[e_{2},y_{3}]=y_{4}$ & \textemdash \\

\hyperref[case:1-4-9]{$F_{1,4}^{(9)}$} & \textemdash & $[e_{2},y_{1}]=y_{2}+y_{3}+2y_{4}$\newline $[e_{2},y_{2}]=y_{3}+y_{4}$\newline $[e_{2},y_{3}]=y_{4}$ & \textemdash \\

\end{longtable}

\begin{longtable}{@{}c>{\raggedright\arraybackslash}p{42mm}>{\raggedright\arraybackslash}p{42mm}>{\raggedright\arraybackslash}p{44mm}@{}}

\caption{Multiplication tables for the family $F_{2,2}$ (superdimension $(3\mid 2)$).}\label{tab:mult-2-2}\\
\toprule
\multicolumn{4}{@{}l@{}}{$[e_{1},e_{2}]=e_{3},\quad [e_{1},y_{1}]=y_{2}$}\\[3pt]
\midrule
Case & $[F_{\bar0},F_{\bar0}]$ & $[F_{\bar0},F_{\bar1}]$ & $[F_{\bar1},F_{\bar1}]$\\
\midrule
\endfirsthead

\multicolumn{4}{@{}l@{}}{\textit{Table \thetable\ (continued; family $F_{2,2}$)}}\\[3pt]
\toprule
\multicolumn{4}{@{}l@{}}{$[e_{1},e_{2}]=e_{3},\quad [e_{1},y_{1}]=y_{2}$}\\[3pt]
\midrule
Case & $[F_{\bar0},F_{\bar0}]$ & $[F_{\bar0},F_{\bar1}]$ & $[F_{\bar1},F_{\bar1}]$\\
\midrule
\endhead

\midrule
\multicolumn{4}{r@{}}{\textit{Continued on the next page}}\\
\endfoot
\bottomrule
\endlastfoot

\hyperref[case:2-2-1]{$F_{2,2}^{(1)}$} & \textemdash & \textemdash & \textemdash \\

\hyperref[case:2-2-2]{$F_{2,2}^{(2)}$} & \textemdash & \textemdash & $[y_{1},y_{1}]=e_{2}$\newline $[y_{1},y_{2}]=\tfrac{1}{2}e_{3}$ \\

\hyperref[case:2-2-3]{$F_{2,2}^{(3)}$} & \textemdash & \textemdash & $[y_{1},y_{1}]=e_{3}$ \\

\hyperref[case:2-2-4]{$F_{2,2}^{(4)}$} & \textemdash & $[e_{2},y_{1}]=y_{2}$ & \textemdash \\

\hyperref[case:2-2-5]{$F_{2,2}^{(5)}$} & \textemdash & $[e_{2},y_{1}]=y_{2}$ & $[y_{1},y_{1}]=e_{3}$ \\

\end{longtable}

\begin{longtable}{@{}c>{\raggedright\arraybackslash}p{42mm}>{\raggedright\arraybackslash}p{42mm}>{\raggedright\arraybackslash}p{44mm}@{}}

\caption{Multiplication tables for the family $F_{2,3}$ (superdimension $(3\mid 3)$).}\label{tab:mult-2-3}\\
\toprule
\multicolumn{4}{@{}l@{}}{$[e_{1},e_{2}]=e_{3},\quad [e_{1},y_{1}]=y_{2},\quad [e_{1},y_{2}]=y_{3}$}\\[3pt]
\midrule
Case & $[F_{\bar0},F_{\bar0}]$ & $[F_{\bar0},F_{\bar1}]$ & $[F_{\bar1},F_{\bar1}]$\\
\midrule
\endfirsthead

\multicolumn{4}{@{}l@{}}{\textit{Table \thetable\ (continued; family $F_{2,3}$)}}\\[3pt]
\toprule
\multicolumn{4}{@{}l@{}}{$[e_{1},e_{2}]=e_{3},\quad [e_{1},y_{1}]=y_{2},\quad [e_{1},y_{2}]=y_{3}$}\\[3pt]
\midrule
Case & $[F_{\bar0},F_{\bar0}]$ & $[F_{\bar0},F_{\bar1}]$ & $[F_{\bar1},F_{\bar1}]$\\
\midrule
\endhead

\midrule
\multicolumn{4}{r@{}}{\textit{Continued on the next page}}\\
\endfoot
\bottomrule
\endlastfoot

\hyperref[case:2-3-1]{$F_{2,3}^{(1)}$} & \textemdash & \textemdash & \textemdash \\

\hyperref[case:2-3-2]{$F_{2,3}^{(2)}$} & \textemdash & $[e_{2},y_{1}]=y_{2}$\newline $[e_{2},y_{2}]=y_{3}$ & \textemdash \\

\hyperref[case:2-3-3]{$F_{2,3}^{(3)}$} & \textemdash & $[e_{2},y_{1}]=y_{3}$ & \textemdash \\

\hyperref[case:2-3-4]{$F_{2,3}^{(4)}$} & \textemdash & $[e_{2},y_{2}]=-y_{3}$\newline $[e_{3},y_{1}]=y_{3}$ & \textemdash \\

\hyperref[case:2-3-5]{$F_{2,3}^{(5)}$} & \textemdash & $[e_{2},y_{1}]=y_{2}+y_{3}$\newline $[e_{2},y_{2}]=y_{3}$ & \textemdash \\

\hyperref[case:2-3-6]{$F_{2,3}^{(6)}$} & \textemdash & $[e_{2},y_{1}]=y_{2}$\newline $[e_{3},y_{1}]=y_{3}$ & \textemdash \\

\hyperref[case:2-3-7]{$F_{2,3}^{(7)}$} & \textemdash & $[e_{2},y_{1}]=2y_{2}$\newline $[e_{2},y_{2}]=y_{3}$\newline $[e_{3},y_{1}]=y_{3}$ & \textemdash \\

\hyperref[case:2-3-8]{$F_{2,3}^{(8)}$} & \textemdash & \textemdash & $[y_{1},y_{1}]=e_{3}$ \\

\hyperref[case:2-3-9]{$F_{2,3}^{(9)}$} & \textemdash & $[e_{2},y_{1}]=y_{2}$\newline $[e_{2},y_{2}]=y_{3}$ & $[y_{1},y_{1}]=e_{3}$ \\

\hyperref[case:2-3-10]{$F_{2,3}^{(10)}$} & \textemdash & $[e_{2},y_{1}]=y_{3}$ & $[y_{1},y_{1}]=e_{3}$ \\

\hyperref[case:2-3-11]{$F_{2,3}^{(11)}$} & \textemdash & $[e_{2},y_{1}]=y_{2}+y_{3}$\newline $[e_{2},y_{2}]=y_{3}$ & $[y_{1},y_{1}]=e_{3}$ \\

\hyperref[case:2-3-12]{$F_{2,3}^{(12)}$} & \textemdash & \textemdash & $[y_{1},y_{3}]=-e_{3}$\newline $[y_{2},y_{2}]=e_{3}$ \\

\hyperref[case:2-3-13]{$F_{2,3}^{(13)}$} & \textemdash & $[e_{2},y_{1}]=y_{2}$\newline $[e_{2},y_{2}]=y_{3}$ & $[y_{1},y_{3}]=-e_{3}$\newline $[y_{2},y_{2}]=e_{3}$ \\

\hyperref[case:2-3-14]{$F_{2,3}^{(14)}$} & \textemdash & \textemdash & $[y_{1},y_{1}]=e_{2}$\newline $[y_{1},y_{2}]=\tfrac{1}{2}e_{3}$ \\

\hyperref[case:2-3-15]{$F_{2,3}^{(15)}$} & \textemdash & $[e_{2},y_{2}]=-y_{3}$\newline $[e_{3},y_{1}]=y_{3}$ & $[y_{1},y_{1}]=e_{2}$\newline $[y_{1},y_{2}]=\tfrac{1}{2}e_{3}$ \\

\hyperref[case:2-3-16]{$F_{2,3}^{(16)}$} & \textemdash & \textemdash & $[y_{1},y_{1}]=e_{2}$\newline $[y_{1},y_{2}]=\tfrac{1}{2}e_{3}$\newline $[y_{1},y_{3}]=-e_{3}$\newline $[y_{2},y_{2}]=e_{3}$ \\

\end{longtable}

\begin{longtable}{@{}c>{\raggedright\arraybackslash}p{42mm}>{\raggedright\arraybackslash}p{42mm}>{\raggedright\arraybackslash}p{44mm}@{}}

\caption{Multiplication tables for the family $F_{3,2}$ (superdimension $(4\mid 2)$).}\label{tab:mult-3-2}\\
\toprule
\multicolumn{4}{@{}l@{}}{$[e_{1},e_{2}]=e_{3},\quad [e_{1},e_{3}]=e_{4},\quad [e_{1},y_{1}]=y_{2}$}\\[3pt]
\midrule
Case & $[F_{\bar0},F_{\bar0}]$ & $[F_{\bar0},F_{\bar1}]$ & $[F_{\bar1},F_{\bar1}]$\\
\midrule
\endfirsthead

\multicolumn{4}{@{}l@{}}{\textit{Table \thetable\ (continued; family $F_{3,2}$)}}\\[3pt]
\toprule
\multicolumn{4}{@{}l@{}}{$[e_{1},e_{2}]=e_{3},\quad [e_{1},e_{3}]=e_{4},\quad [e_{1},y_{1}]=y_{2}$}\\[3pt]
\midrule
Case & $[F_{\bar0},F_{\bar0}]$ & $[F_{\bar0},F_{\bar1}]$ & $[F_{\bar1},F_{\bar1}]$\\
\midrule
\endhead

\midrule
\multicolumn{4}{r@{}}{\textit{Continued on the next page}}\\
\endfoot
\bottomrule
\endlastfoot

\hyperref[case:3-2-1]{$F_{3,2}^{(1)}$} & \textemdash & \textemdash & \textemdash \\

\hyperref[case:3-2-2]{$F_{3,2}^{(2)}$} & \textemdash & \textemdash & $[y_{1},y_{1}]=e_{2}$\newline $[y_{1},y_{2}]=\tfrac{1}{2}e_{3}$\newline $[y_{2},y_{2}]=\tfrac{1}{2}e_{4}$ \\

\hyperref[case:3-2-3]{$F_{3,2}^{(3)}$} & \textemdash & \textemdash & $[y_{1},y_{1}]=e_{3}$\newline $[y_{1},y_{2}]=\tfrac{1}{2}e_{4}$ \\

\hyperref[case:3-2-4]{$F_{3,2}^{(4)}$} & \textemdash & \textemdash & $[y_{1},y_{1}]=e_{4}$ \\

\hyperref[case:3-2-5]{$F_{3,2}^{(5)}$} & \textemdash & $[e_{2},y_{1}]=y_{2}$ & \textemdash \\

\hyperref[case:3-2-6]{$F_{3,2}^{(6)}$} & \textemdash & $[e_{2},y_{1}]=y_{2}$ & $[y_{1},y_{1}]=e_{4}$ \\

\end{longtable}

\begin{longtable}{@{}c>{\raggedright\arraybackslash}p{42mm}>{\raggedright\arraybackslash}p{42mm}>{\raggedright\arraybackslash}p{44mm}@{}}

\caption{Multiplication tables for the family $F_{4,2}$ (superdimension $(5\mid 2)$).}\label{tab:mult-4-2}\\
\toprule
\multicolumn{4}{@{}l@{}}{$[e_{1},e_{2}]=e_{3},\quad [e_{1},e_{3}]=e_{4},\quad [e_{1},e_{4}]=e_{5},\quad [e_{1},y_{1}]=y_{2}$}\\[3pt]
\midrule
Case & $[F_{\bar0},F_{\bar0}]$ & $[F_{\bar0},F_{\bar1}]$ & $[F_{\bar1},F_{\bar1}]$\\
\midrule
\endfirsthead

\multicolumn{4}{@{}l@{}}{\textit{Table \thetable\ (continued; family $F_{4,2}$)}}\\[3pt]
\toprule
\multicolumn{4}{@{}l@{}}{$[e_{1},e_{2}]=e_{3},\quad [e_{1},e_{3}]=e_{4},\quad [e_{1},e_{4}]=e_{5},\quad [e_{1},y_{1}]=y_{2}$}\\[3pt]
\midrule
Case & $[F_{\bar0},F_{\bar0}]$ & $[F_{\bar0},F_{\bar1}]$ & $[F_{\bar1},F_{\bar1}]$\\
\midrule
\endhead

\midrule
\multicolumn{4}{r@{}}{\textit{Continued on the next page}}\\
\endfoot
\bottomrule
\endlastfoot

\hyperref[case:4-2-1]{$F_{4,2}^{(1)}$} & \textemdash & \textemdash & \textemdash \\

\hyperref[case:4-2-2]{$F_{4,2}^{(2)}$} & \textemdash & \textemdash & $[y_{1},y_{1}]=e_{3}$\newline $[y_{1},y_{2}]=\tfrac{1}{2}e_{4}$\newline $[y_{2},y_{2}]=\tfrac{1}{2}e_{5}$ \\

\hyperref[case:4-2-3]{$F_{4,2}^{(3)}$} & \textemdash & \textemdash & $[y_{1},y_{1}]=e_{4}$\newline $[y_{1},y_{2}]=\tfrac{1}{2}e_{5}$ \\

\hyperref[case:4-2-4]{$F_{4,2}^{(4)}$} & \textemdash & \textemdash & $[y_{1},y_{1}]=e_{5}$ \\

\hyperref[case:4-2-5]{$F_{4,2}^{(5)}$} & \textemdash & $[e_{2},y_{1}]=y_{2}$ & \textemdash \\

\hyperref[case:4-2-6]{$F_{4,2}^{(6)}$} & \textemdash & $[e_{2},y_{1}]=y_{2}$ & $[y_{1},y_{1}]=e_{5}$ \\

\hyperref[case:4-2-7]{$F_{4,2}^{(7)}$} & $[e_{2},e_{3}]=e_{5}$ & \textemdash & \textemdash \\

\hyperref[case:4-2-8]{$F_{4,2}^{(8)}$} & $[e_{2},e_{3}]=e_{5}$ & \textemdash & $[y_{1},y_{1}]=e_{4}$\newline $[y_{1},y_{2}]=\tfrac{1}{2}e_{5}$ \\

\hyperref[case:4-2-9]{$F_{4,2}^{(9)}$} & $[e_{2},e_{3}]=e_{5}$ & \textemdash & $[y_{1},y_{1}]=e_{5}$ \\

\hyperref[case:4-2-10]{$F_{4,2}^{(10)}$} & $[e_{2},e_{3}]=e_{5}$ & $[e_{2},y_{1}]=y_{2}$ & \textemdash \\

\hyperref[case:4-2-11]{$F_{4,2}^{(11)}$} & $[e_{2},e_{3}]=e_{5}$ & $[e_{2},y_{1}]=y_{2}$ & $[y_{1},y_{1}]=e_{5}$ \\

\end{longtable}

\begin{longtable}{@{}c>{\raggedright\arraybackslash}p{42mm}>{\raggedright\arraybackslash}p{42mm}>{\raggedright\arraybackslash}p{44mm}@{}}

\caption{Multiplication tables for the family $F_{5,2}$ (superdimension $(6\mid 2)$).}\label{tab:mult-5-2}\\
\toprule
\multicolumn{4}{@{}l@{}}{$[e_{1},e_{2}]=e_{3},\quad [e_{1},e_{3}]=e_{4},\quad [e_{1},e_{4}]=e_{5},\quad [e_{1},e_{5}]=e_{6},\quad [e_{1},y_{1}]=y_{2}$}\\[3pt]
\midrule
Case & $[F_{\bar0},F_{\bar0}]$ & $[F_{\bar0},F_{\bar1}]$ & $[F_{\bar1},F_{\bar1}]$\\
\midrule
\endfirsthead

\multicolumn{4}{@{}l@{}}{\textit{Table \thetable\ (continued; family $F_{5,2}$)}}\\[3pt]
\toprule
\multicolumn{4}{@{}l@{}}{$[e_{1},e_{2}]=e_{3},\quad [e_{1},e_{3}]=e_{4},\quad [e_{1},e_{4}]=e_{5},\quad [e_{1},e_{5}]=e_{6},\quad [e_{1},y_{1}]=y_{2}$}\\[3pt]
\midrule
Case & $[F_{\bar0},F_{\bar0}]$ & $[F_{\bar0},F_{\bar1}]$ & $[F_{\bar1},F_{\bar1}]$\\
\midrule
\endhead

\midrule
\multicolumn{4}{r@{}}{\textit{Continued on the next page}}\\
\endfoot
\bottomrule
\endlastfoot

\hyperref[case:5-2-1]{$F_{5,2}^{(1)}$} & \textemdash & \textemdash & \textemdash \\

\hyperref[case:5-2-2]{$F_{5,2}^{(2)}$} & \textemdash & \textemdash & $[y_{1},y_{1}]=e_{4}$\newline $[y_{1},y_{2}]=\tfrac{1}{2}e_{5}$\newline $[y_{2},y_{2}]=\tfrac{1}{2}e_{6}$ \\

\hyperref[case:5-2-3]{$F_{5,2}^{(3)}$} & \textemdash & \textemdash & $[y_{1},y_{1}]=e_{5}$\newline $[y_{1},y_{2}]=\tfrac{1}{2}e_{6}$ \\

\hyperref[case:5-2-4]{$F_{5,2}^{(4)}$} & \textemdash & \textemdash & $[y_{1},y_{1}]=e_{6}$ \\

\hyperref[case:5-2-5]{$F_{5,2}^{(5)}$} & \textemdash & $[e_{2},y_{1}]=y_{2}$ & \textemdash \\

\hyperref[case:5-2-6]{$F_{5,2}^{(6)}$} & \textemdash & $[e_{2},y_{1}]=y_{2}$ & $[y_{1},y_{1}]=e_{6}$ \\

\hyperref[case:5-2-7]{$F_{5,2}^{(7)}$} & $[e_{2},e_{3}]=e_{6}$ & \textemdash & \textemdash \\

\hyperref[case:5-2-8]{$F_{5,2}^{(8)}$} & $[e_{2},e_{3}]=e_{6}$ & \textemdash & $[y_{1},y_{1}]=e_{4}$\newline $[y_{1},y_{2}]=\tfrac{1}{2}e_{5}$\newline $[y_{2},y_{2}]=\tfrac{1}{2}e_{6}$ \\

\hyperref[case:5-2-9]{$F_{5,2}^{(9)}$} & $[e_{2},e_{3}]=e_{6}$ & \textemdash & $[y_{1},y_{1}]=e_{5}$\newline $[y_{1},y_{2}]=\tfrac{1}{2}e_{6}$ \\

\hyperref[case:5-2-10]{$F_{5,2}^{(10)}$} & $[e_{2},e_{3}]=e_{6}$ & \textemdash & $[y_{1},y_{1}]=e_{6}$ \\

\hyperref[case:5-2-11]{$F_{5,2}^{(11)}$} & $[e_{2},e_{3}]=e_{6}$ & $[e_{2},y_{1}]=y_{2}$ & \textemdash \\

\hyperref[case:5-2-12]{$F_{5,2}^{(12)}$} & $[e_{2},e_{3}]=e_{6}$ & $[e_{2},y_{1}]=y_{2}$ & $[y_{1},y_{1}]=e_{6}$ \\

\hyperref[case:5-2-13]{$F_{5,2}^{(13)}$} & $[e_{2},e_{3}]=e_{5}$\newline $[e_{2},e_{4}]=e_{6}$\newline $[e_{2},e_{5}]=-e_{6}$\newline $[e_{3},e_{4}]=e_{6}$ & \textemdash & \textemdash \\

\hyperref[case:5-2-14]{$F_{5,2}^{(14)}$} & $[e_{2},e_{3}]=e_{5}$\newline $[e_{2},e_{4}]=e_{6}$\newline $[e_{2},e_{5}]=-e_{6}$\newline $[e_{3},e_{4}]=e_{6}$ & \textemdash & $[y_{1},y_{1}]=e_{6}$ \\

\hyperref[case:5-2-15]{$F_{5,2}^{(15)}$} & $[e_{2},e_{3}]=e_{5}$\newline $[e_{2},e_{4}]=e_{6}$\newline $[e_{2},e_{5}]=-e_{6}$\newline $[e_{3},e_{4}]=e_{6}$ & $[e_{2},y_{1}]=y_{2}$ & \textemdash \\

\hyperref[case:5-2-16]{$F_{5,2}^{(16)}$} & $[e_{2},e_{3}]=e_{5}$\newline $[e_{2},e_{4}]=e_{6}$\newline $[e_{2},e_{5}]=-e_{6}$\newline $[e_{3},e_{4}]=e_{6}$ & $[e_{2},y_{1}]=y_{2}$ & $[y_{1},y_{1}]=e_{6}$ \\

\hyperref[case:5-2-17]{$F_{5,2}^{(17)}$} & $[e_{2},e_{5}]=-e_{6}$\newline $[e_{3},e_{4}]=e_{6}$ & \textemdash & \textemdash \\

\hyperref[case:5-2-18]{$F_{5,2}^{(18)}$} & $[e_{2},e_{5}]=-e_{6}$\newline $[e_{3},e_{4}]=e_{6}$ & \textemdash & $[y_{1},y_{1}]=e_{6}$ \\

\hyperref[case:5-2-19]{$F_{5,2}^{(19)}$} & $[e_{2},e_{5}]=-e_{6}$\newline $[e_{3},e_{4}]=e_{6}$ & $[e_{2},y_{1}]=y_{2}$ & \textemdash \\

\hyperref[case:5-2-20]{$F_{5,2}^{(20)}$} & $[e_{2},e_{5}]=-e_{6}$\newline $[e_{3},e_{4}]=e_{6}$ & $[e_{2},y_{1}]=y_{2}$ & $[y_{1},y_{1}]=e_{6}$ \\

\hyperref[case:5-2-21]{$F_{5,2}^{(21)}$} & $[e_{2},e_{5}]=-e_{6}$\newline $[e_{3},e_{4}]=e_{6}$ & $[e_{2},y_{1}]=-y_{2}$ & \textemdash \\

\hyperref[case:5-2-22]{$F_{5,2}^{(22)}$} & $[e_{2},e_{5}]=-e_{6}$\newline $[e_{3},e_{4}]=e_{6}$ & $[e_{2},y_{1}]=-y_{2}$ & $[y_{1},y_{1}]=e_{5}$\newline $[y_{1},y_{2}]=\tfrac{1}{2}e_{6}$ \\

\hyperref[case:5-2-23]{$F_{5,2}^{(23)}$} & $[e_{2},e_{5}]=-e_{6}$\newline $[e_{3},e_{4}]=e_{6}$ & $[e_{2},y_{1}]=-y_{2}$ & $[y_{1},y_{1}]=e_{6}$ \\

\hyperref[case:5-2-24]{$F_{5,2}^{(24)}$} & $[e_{2},e_{3}]=e_{6}$\newline $[e_{2},e_{5}]=-e_{6}$\newline $[e_{3},e_{4}]=e_{6}$ & \textemdash & \textemdash \\

\hyperref[case:5-2-25]{$F_{5,2}^{(25)}$} & $[e_{2},e_{3}]=e_{6}$\newline $[e_{2},e_{5}]=-e_{6}$\newline $[e_{3},e_{4}]=e_{6}$ & \textemdash & $[y_{1},y_{1}]=e_{6}$ \\

\hyperref[case:5-2-26]{$F_{5,2}^{(26)}$} & $[e_{2},e_{3}]=e_{6}$\newline $[e_{2},e_{5}]=-e_{6}$\newline $[e_{3},e_{4}]=e_{6}$ & $[e_{2},y_{1}]=y_{2}$ & \textemdash \\

\hyperref[case:5-2-27]{$F_{5,2}^{(27)}$} & $[e_{2},e_{3}]=e_{6}$\newline $[e_{2},e_{5}]=-e_{6}$\newline $[e_{3},e_{4}]=e_{6}$ & $[e_{2},y_{1}]=y_{2}$ & $[y_{1},y_{1}]=e_{6}$ \\

\hyperref[case:5-2-28]{$F_{5,2}^{(28)}$} & $[e_{2},e_{3}]=e_{6}$\newline $[e_{2},e_{5}]=-e_{6}$\newline $[e_{3},e_{4}]=e_{6}$ & $[e_{2},y_{1}]=-y_{2}$ & \textemdash \\

\hyperref[case:5-2-29]{$F_{5,2}^{(29)}$} & $[e_{2},e_{3}]=e_{6}$\newline $[e_{2},e_{5}]=-e_{6}$\newline $[e_{3},e_{4}]=e_{6}$ & $[e_{2},y_{1}]=-y_{2}$ & $[y_{1},y_{1}]=e_{5}$\newline $[y_{1},y_{2}]=\tfrac{1}{2}e_{6}$ \\

\hyperref[case:5-2-30]{$F_{5,2}^{(30)}$} & $[e_{2},e_{3}]=e_{6}$\newline $[e_{2},e_{5}]=-e_{6}$\newline $[e_{3},e_{4}]=e_{6}$ & $[e_{2},y_{1}]=-y_{2}$ & $[y_{1},y_{1}]=e_{6}$ \\

\hyperref[case:5-2-31]{$F_{5,2}^{(31)}$} & $[e_{2},e_{3}]=e_{5}+e_{6}$\newline $[e_{2},e_{4}]=e_{6}$ & $[e_{2},y_{1}]=-2y_{2}$ & \textemdash \\

\hyperref[case:5-2-32]{$F_{5,2}^{(32)}$} & $[e_{2},e_{3}]=e_{5}+e_{6}$\newline $[e_{2},e_{4}]=e_{6}$ & $[e_{2},y_{1}]=-2y_{2}$ & $[y_{1},y_{1}]=e_{6}$ \\

\hyperref[case:5-2-33]{$F_{5,2}^{(33)}$} & $[e_{2},e_{3}]=e_{5}$\newline $[e_{2},e_{4}]=e_{6}$\newline $[e_{2},e_{5}]=-e_{6}$\newline $[e_{3},e_{4}]=e_{6}$ & $[e_{2},y_{1}]=-y_{2}$ & \textemdash \\

\hyperref[case:5-2-34]{$F_{5,2}^{(34)}$} & $[e_{2},e_{3}]=e_{5}$\newline $[e_{2},e_{4}]=e_{6}$\newline $[e_{2},e_{5}]=-e_{6}$\newline $[e_{3},e_{4}]=e_{6}$ & $[e_{2},y_{1}]=-y_{2}$ & $[y_{1},y_{1}]=e_{5}$\newline $[y_{1},y_{2}]=\tfrac{1}{2}e_{6}$ \\

\hyperref[case:5-2-35]{$F_{5,2}^{(35)}$} & $[e_{2},e_{3}]=e_{5}$\newline $[e_{2},e_{4}]=e_{6}$\newline $[e_{2},e_{5}]=-e_{6}$\newline $[e_{3},e_{4}]=e_{6}$ & $[e_{2},y_{1}]=-y_{2}$ & $[y_{1},y_{1}]=e_{6}$ \\

\end{longtable}

\endgroup

\subsection{Homogeneous derivations}

A homogeneous linear map $d$ of parity $\varepsilon\in\{0,1\}$ is called a derivation if $d(F_{\bar\alpha})\subseteq F_{\overline{\alpha+\varepsilon}}$ and, for all homogeneous $x,y\in F$,
\begin{equation}\label{eq:der}
 d([x,y])=[d(x),y]+(-1)^{\varepsilon|x|}[x,d(y)].
\end{equation}
Write the ordered basis \eqref{eq:basis} as $(x_1,\ldots,x_{p+q})$. We adopt the row convention:
\begin{equation}\label{eq:row}
 d(x_i)=\sum_{j=1}^{p+q}D_{ij}x_j.
\end{equation}
Thus, row $i$ contains the coordinates of the image $d(x_i)$. In the standard column convention, the corresponding matrix is $D^{\mathsf T}$. In particular, even and odd derivations have the respective forms
\begin{equation}\label{eq:blocks}
 D_{\bar0}=\begin{pmatrix} A&0\\0&B\end{pmatrix},\qquad
 D_{\bar1}=\begin{pmatrix}0&U\\V&0\end{pmatrix},
\end{equation}
where $A\in\mathbb C^{p\times p}$, $B\in\mathbb C^{q\times q}$, $U\in\mathbb C^{p\times q}$ and $V\in\mathbb C^{q\times p}$.

\section{Computational method}
\begin{lemma}\label{lem:kernel}
Let $[x_i,x_j]=\sum_k C_{ij}^k x_k$ and $N=p+q$. For $\varepsilon=0,1$, the space $\Der_{\bar\varepsilon}(F)$ is isomorphic to the kernel of the homogeneous linear system
\begin{equation}\label{eq:linear}
 \sum_{\ell=1}^{N}\left(C_{ij}^{\ell}D_{\ell k}-D_{i\ell}C_{\ell j}^{k}
 -(-1)^{\varepsilon |x_i|}D_{j\ell}C_{i\ell}^{k}\right)=0,
\end{equation}
where $1\le i\le j\le N$, $1\le k\le N$, and $D_{ij}=0$ whenever $|x_j|\not\equiv |x_i|+\varepsilon\pmod 2$. If $L_\varepsilon$ denotes the coefficient matrix of this system after the entries excluded by parity have been removed, then
\begin{equation}\label{eq:rank}
\begin{aligned}
 \dim\Der_{\bar0}(F)&=p^2+q^2-\rk L_0,\\
 \dim\Der_{\bar1}(F)&=2pq-\rk L_1.
\end{aligned}
\end{equation}
\end{lemma}
\begin{proof}
Substituting \eqref{eq:row} and the structure constants into \eqref{eq:der} and comparing the coefficients of $x_k$ gives \eqref{eq:linear}. Conversely, if the system is satisfied, the identity \eqref{eq:der} holds for every ordered pair of basis elements: the pairs with $i>j$ follow by super-skew-symmetry. Bilinearity on each pair of homogeneous components extends the identity to all homogeneous elements. Hence, the system is necessary and sufficient. The formulas in \eqref{eq:rank} follow from the rank--nullity theorem.
\end{proof}

For a direct computation, it suffices to specify the images of $e_1,e_2,y_1$, since the remaining images are determined recursively by \eqref{eq:chains}:
\begin{equation}\label{eq:rec}
\begin{aligned}
 d(e_{i+1})&=[d(e_1),e_i]+[e_1,d(e_i)],&& 2\le i\le n,\\
 d(y_{j+1})&=[d(e_1),y_j]+[e_1,d(y_j)],&& 1\le j<m.
\end{aligned}
\end{equation}
The derivation identity is then checked for the remaining basis brackets, including the zero brackets. For an even derivation $d$, the coordinates of the images of $e_1,e_2,y_1$ are denoted by $a_i,b_i,c_j$, respectively; for an odd derivation, they are denoted by $u_j,v_j,w_i$, respectively.

\begin{example}\label{ex:f123}
Consider $F=F_{1,2}^{(3)}$. Its nonzero brackets, up to super-skew-symmetry, are $[e_1,y_1]=y_2$ and $[y_1,y_1]=e_2$.

\emph{Even case.} Set
\[
\begin{aligned}
 d(e_1)&=a_1e_1+a_2e_2,&d(e_2)&=b_1e_1+b_2e_2,\\
 d(y_1)&=c_1y_1+c_2y_2,&d(y_2)&=h_1y_1+h_2y_2.
\end{aligned}
\]
Applying \eqref{eq:der} to the two nonzero brackets gives
\[
\begin{aligned}
 d(y_2)&=[a_1e_1+a_2e_2,y_1]+[e_1,c_1y_1+c_2y_2]=(a_1+c_1)y_2,\\
 d(e_2)&=2[d(y_1),y_1]=2c_1e_2.
\end{aligned}
\]
Consequently, $h_1=b_1=0$, $h_2=a_1+c_1$ and $b_2=2c_1$.

\emph{Odd case.} Let
\[
\begin{aligned}
 d(e_1)&=u_1y_1+u_2y_2,&d(e_2)&=v_1y_1+v_2y_2,\\
 d(y_1)&=w_1e_1+w_2e_2,&d(y_2)&=z_1e_1+z_2e_2.
\end{aligned}
\]
Then
\[
\begin{aligned}
 d(y_2)&=[u_1y_1+u_2y_2,y_1]+[e_1,w_1e_1+w_2e_2]=u_1e_2,\\
 d(e_2)&=[d(y_1),y_1]-[y_1,d(y_1)]=2w_1y_2.
\end{aligned}
\]
Hence, $z_1=v_1=0$, $z_2=u_1$ and $v_2=2w_1$.
In both cases, the images of $e_2,y_2$ lie in the center. For the remaining pairs $(e_1,e_2)$, $(e_1,y_2)$, $(e_2,y_1)$, $(e_2,y_2)$, $(y_1,y_2)$ and $(y_2,y_2)$, both sides of \eqref{eq:der} vanish; for equal even arguments, the identity holds automatically. Thus, there are no further restrictions, and we obtain
\[
D_{\bar0}=\begin{pmatrix}a_1&a_2&0&0\\0&2c_1&0&0\\0&0&c_1&c_2\\0&0&0&a_1+c_1\end{pmatrix},\qquad
D_{\bar1}=\begin{pmatrix}0&0&u_1&u_2\\0&0&0&2w_1\\w_1&w_2&0&0\\0&u_1&0&0\end{pmatrix}.
\]
Each case has four independent parameters, so $\dim\Der_{\bar0}(F)=\dim\Der_{\bar1}(F)=4$.
\end{example}

\section{Main results}
\begin{theorem}\label{thm:main}
For each $F_{n,m}^{(i)}$ in the eight families under consideration, all even and odd derivations are given by the matrices \eqref{eq:blocks}, whose blocks are listed in Appendix~\ref{app:matrices}. For each parity and each case, the parameters appearing in the corresponding pair of blocks are arbitrary and independent over $\mathbb C$. The dimensions of the homogeneous components are given in Table~\ref{tab:dims}.
\end{theorem}
\begin{proof}
For each multiplication table in Section~2.1, we construct two systems of the form \eqref{eq:linear}. All structure constants, including the zero ones, are taken into account; the super-Jacobi identity is also verified for the original tables. Since the structure constants are rational, the reduced row echelon form of each system is computed over $\mathbb Q$ using exact arithmetic. Extending the field to $\mathbb C$ does not change the rank. The computations can be reproduced using the accompanying Python script (SymPy).

To obtain the parametrizations, we first eliminate the coordinates of the images of basis elements other than $e_1,e_2,y_1$, and then the dependent coordinates of these three images. The remaining free coordinates are denoted as in Section~3. Substituting the resulting matrices into \eqref{eq:linear} yields zero on the left-hand side for every pair and every coefficient. The number of free coordinates equals the nullity of the corresponding system. By Lemma~\ref{lem:kernel}, this establishes both that the indicated maps are derivations and that the parametrizations are complete. The results of the elimination are given in the appendix, and the numbers of free coordinates are listed in Table~\ref{tab:dims}.
\end{proof}
\begin{table}[htbp]
\centering\small
\caption{Dimensions $d_0(i)=\dim\Der_{\bar0}(F_{n,m}^{(i)})$ and $d_1(i)=\dim\Der_{\bar1}(F_{n,m}^{(i)})$.}
\label{tab:dims}
\setlength{\tabcolsep}{5pt}
\begin{tabular}{@{}clll@{}}
\toprule
Family & Indices $i$ & $d_0(i)$ & $d_1(i)$\\
\midrule
$F_{1,2}$ & 1--3 & $5,5,4$ & $4,4,4$ \\
\addlinespace[2pt]
$F_{1,3}$ & 1--6 & $6,5,4,6,5,5$ & $5,4,3,5,4,4$ \\
\addlinespace[2pt]
$F_{1,4}$ & 1--9 & $7,6,5,7,5,6,5,6,5$ & $6,5,5,6,5,5,5,5,5$ \\
\addlinespace[2pt]
$F_{2,2}$ & 1--5 & $7,5,6,7,6$ & $6,4,5,6,5$ \\
\addlinespace[2pt]
$F_{2,3}$ & 1--9 & $8,8,7,6,7,6,6,7,7$ & $7,7,6,5,6,5,5,6,6$ \\
 & 10--16 & $6,6,6,6,6,5,5$ & $6,6,4,4,6,5,4$ \\
\addlinespace[2pt]
$F_{3,2}$ & 1--6 & $9,6,7,8,8,7$ & $6,3,4,5,5,5$ \\
\addlinespace[2pt]
$F_{4,2}$ & 1--9 & $11,8,9,10,10,9,10,8,9$ & $6,4,4,5,5,5,6,4,5$ \\
 & 10--11 & $9,8$ & $5,5$ \\
\addlinespace[2pt]
$F_{5,2}$ & 1--9 & $13,10,11,12,12,11,12,9,10$ & $6,4,4,5,5,5,6,3,4$ \\
 & 10--18 & $11,11,10,10,9,9,8,11,10$ & $5,5,5,5,5,5,5,5,5$ \\
 & 19--27 & $10,9,11,9,10,11,10,10,9$ & $5,5,6,3,5,5,5,5,5$ \\
 & 28--35 & $11,9,10,11,10,10,8,9$ & $6,3,5,5,5,6,3,5$ \\
\addlinespace[2pt]
\bottomrule
\end{tabular}
\end{table}

\begin{corollary}
For all 91 cases under consideration, the following inequalities hold:
\[
4\le\dim\Der_{\bar0}(F)\le13,\qquad
3\le\dim\Der_{\bar1}(F)\le7,\qquad
7\le\dim\Der(F)\le19.
\]
The minimum total dimension of $7$ is attained by $F_{1,3}^{(3)}$, and the maximum of $19$ by $F_{5,2}^{(1)}$. In particular, the odd component is nonzero in every case under consideration.
\end{corollary}
\begin{proof}
The assertion follows directly by comparing the dimensions and their sums in Table~\ref{tab:dims}.
\end{proof}

\section{Conclusion}
We have obtained an explicit description of even and odd derivations for all cases in the eight selected families of complex filiform Lie superalgebras. The unified linear system makes the computations verifiable, while the block-matrix representation separates the essential coefficients from the zero blocks imposed by the grading. These results can be used in further investigations of solvable extensions with a filiform nilradical; the classification of the extensions themselves is beyond the scope of this paper.

\par\medskip
\noindent\textbf{Computational and language assistance.} Python with the SymPy library was used for symbolic computations and for checking intermediate calculations and matrix parametrizations. ChatGPT (OpenAI) was employed as an auxiliary tool for computational assistance, LaTeX preparation, English translation, and linguistic editing of the manuscript. The authors retain responsibility for the mathematical content, including all proofs, computations, and results, and for its final verification.

\medskip
\noindent\textbf{Conflict of Interest.} The authors declare that there is no conflict of interest regarding the publication of this article.

\Needspace{9\baselineskip}

\newpage

\newsavebox{\allblocksbox}
\newsavebox{\firstblocksbox}
\newcommand{\derivcase}[6]{%
\par\noindent\begin{minipage}{\linewidth}\small
\setlength{\arraycolsep}{2.5pt}\renewcommand{\arraystretch}{0.98}%
\phantomsection\label{#1}%
\sbox{\allblocksbox}{$\displaystyle #2.\quad #3,\; #4,\; #5,\; #6$}%
\ifdim\wd\allblocksbox<1.111\linewidth
  \typeout{ONE-ROW: #1}%
  \noindent\makebox[\linewidth][l]{\adjustbox{max width=\linewidth}{\usebox{\allblocksbox}}}\par
\else
  \sbox{\allblocksbox}{$\displaystyle #3,\; #4,\; #5,\; #6$}%
  \ifdim\wd\allblocksbox<1.111\linewidth
    \typeout{ONE-ROW: #1}%
    \noindent\textbf{$#2$.}\par\vspace{1pt}%
    \noindent\makebox[\linewidth][c]{\adjustbox{max width=\linewidth}{\usebox{\allblocksbox}}}\par
  \else
    \typeout{TWO-ROW: #1}%
    \sbox{\firstblocksbox}{$\displaystyle #2.\quad #3,\; #5$}%
    \ifdim\wd\firstblocksbox<1.111\linewidth
      \noindent\makebox[\linewidth][l]{\adjustbox{max width=\linewidth}{\usebox{\firstblocksbox}}}\par
    \else
      \noindent\textbf{$#2$.}\par\vspace{1pt}%
      \noindent\makebox[\linewidth][c]{\fitmath{#3,\quad #5}}\par
    \fi
    \vspace{2pt}\noindent\makebox[\linewidth][c]{\fitmath{#4,\quad #6}}\par
  \fi
\fi
\end{minipage}\par\vspace{5pt}%
}
\appendix
\section{Matrix parametrizations}\label{app:matrices}
For each case, we list the four blocks $A,B,U,V$ in \eqref{eq:blocks} with respect to the basis \eqref{eq:basis}. The full matrix of an even derivation is obtained from $A,B$ by inserting two zero off-diagonal blocks, whereas that of an odd derivation is obtained from $U,V$ by inserting zero diagonal blocks. Compact blocks are displayed on a single line; wider matrices are arranged in two lines ($A,U$ and $B,V$). The parameters $a_i,b_i,c_j$ pertain exclusively to even derivations, and $u_j,v_j,w_i$ to odd derivations. Within each case and parity, the parameters are independent; the use of the same symbols in different cases imposes no relations between them.
\subsection{Family \texorpdfstring{$F_{1,2}$}{F1,2} \texorpdfstring{(superdimension $(2\mid2)$)}{(superdimension 2|2)}}
\derivcase{case:1-2-1}{F_{1,2}^{(1)}}
{A=\begin{pmatrix}a_{1} & a_{2} \\ 0 & b_{2}\end{pmatrix}}
{B=\begin{pmatrix}c_{1} & c_{2} \\ 0 & a_{1} + c_{1}\end{pmatrix}}
{U=\begin{pmatrix}u_{1} & u_{2} \\ 0 & v_{2}\end{pmatrix}}
{V=\begin{pmatrix}0 & w_{2} \\ 0 & 0\end{pmatrix}}
\derivcase{case:1-2-2}{F_{1,2}^{(2)}}
{A=\begin{pmatrix}a_{1} & a_{2} \\ b_{1} & a_{1} + a_{2} - b_{1}\end{pmatrix}}
{B=\begin{pmatrix}c_{1} & c_{2} \\ 0 & a_{1} + a_{2} + c_{1}\end{pmatrix}}
{U=\begin{pmatrix}u_{1} & u_{2} \\ u_{1} & v_{2}\end{pmatrix}}
{V=\begin{pmatrix}w_{1} & - w_{1} \\ 0 & 0\end{pmatrix}}
\derivcase{case:1-2-3}{F_{1,2}^{(3)}}
{A=\begin{pmatrix}a_{1} & a_{2} \\ 0 & 2 c_{1}\end{pmatrix}}
{B=\begin{pmatrix}c_{1} & c_{2} \\ 0 & a_{1} + c_{1}\end{pmatrix}}
{U=\begin{pmatrix}u_{1} & u_{2} \\ 0 & 2 w_{1}\end{pmatrix}}
{V=\begin{pmatrix}w_{1} & w_{2} \\ 0 & u_{1}\end{pmatrix}}
\subsection{Family \texorpdfstring{$F_{1,3}$}{F1,3} \texorpdfstring{(superdimension $(2\mid3)$)}{(superdimension 2|3)}}
\derivcase{case:1-3-1}{F_{1,3}^{(1)}}
{A=\begin{pmatrix}a_{1} & a_{2} \\ 0 & b_{2}\end{pmatrix}}
{B=\begin{pmatrix}c_{1} & c_{2} & c_{3} \\ 0 & a_{1} + c_{1} & c_{2} \\ 0 & 0 & 2 a_{1} + c_{1}\end{pmatrix}}
{U=\begin{pmatrix}u_{1} & u_{2} & u_{3} \\ 0 & 0 & v_{3}\end{pmatrix}}
{V=\begin{pmatrix}0 & w_{2} \\ 0 & 0 \\ 0 & 0\end{pmatrix}}
\derivcase{case:1-3-2}{F_{1,3}^{(2)}}
{A=\begin{pmatrix}a_{1} & a_{2} \\ 0 & 2 c_{1}\end{pmatrix}}
{B=\begin{pmatrix}c_{1} & c_{2} & c_{3} \\ 0 & a_{1} + c_{1} & c_{2} \\ 0 & 0 & 2 a_{1} + c_{1}\end{pmatrix}}
{U=\begin{pmatrix}u_{1} & u_{2} & u_{3} \\ 0 & 0 & 0\end{pmatrix}}
{V=\begin{pmatrix}0 & w_{2} \\ 0 & u_{1} \\ 0 & 0\end{pmatrix}}
\derivcase{case:1-3-3}{F_{1,3}^{(3)}}
{A=\begin{pmatrix}a_{1} & a_{2} \\ 0 & 2 \left(a_{1} + c_{1}\right)\end{pmatrix}}
{B=\begin{pmatrix}c_{1} & c_{2} & 0 \\ 0 & a_{1} + c_{1} & c_{2} \\ 0 & 0 & 2 a_{1} + c_{1}\end{pmatrix}}
{U=\begin{pmatrix}0 & u_{2} & u_{3} \\ 0 & 0 & 0\end{pmatrix}}
{V=\begin{pmatrix}0 & w_{2} \\ 0 & - u_{3} \\ 0 & u_{2}\end{pmatrix}}
\derivcase{case:1-3-4}{F_{1,3}^{(4)}}
{A=\begin{pmatrix}a_{1} & a_{2} \\ b_{1} & a_{1} + a_{2} - b_{1}\end{pmatrix}}
{B=\begin{pmatrix}c_{1} & c_{2} & c_{3} \\ 0 & a_{1} + a_{2} + c_{1} & c_{2} \\ 0 & 0 & 2 a_{1} + 2 a_{2} + c_{1}\end{pmatrix}}
{U=\begin{pmatrix}u_{1} & u_{2} & u_{3} \\ u_{1} & u_{2} & v_{3}\end{pmatrix}}
{V=\begin{pmatrix}w_{1} & - w_{1} \\ 0 & 0 \\ 0 & 0\end{pmatrix}}
\derivcase{case:1-3-5}{F_{1,3}^{(5)}}
{A=\begin{pmatrix}a_{1} & a_{2} \\ 0 & 2 a_{1}\end{pmatrix}}
{B=\begin{pmatrix}c_{1} & c_{2} & c_{3} \\ 0 & a_{1} + c_{1} & a_{2} + c_{2} \\ 0 & 0 & 2 a_{1} + c_{1}\end{pmatrix}}
{U=\begin{pmatrix}u_{1} & u_{2} & u_{3} \\ 0 & u_{1} & v_{3}\end{pmatrix}}
{V=\begin{pmatrix}0 & 0 \\ 0 & 0 \\ 0 & 0\end{pmatrix}}
\derivcase{case:1-3-6}{F_{1,3}^{(6)}}
{A=\begin{pmatrix}a_{1} & a_{2} \\ - a_{1} - 2 a_{2} & 2 a_{1} + 3 a_{2}\end{pmatrix}}
{B=\begin{pmatrix}c_{1} & c_{2} & c_{3} \\ 0 & a_{1} + a_{2} + c_{1} & a_{2} + c_{2} \\ 0 & 0 & 2 a_{1} + 2 a_{2} + c_{1}\end{pmatrix}}
{U=\begin{pmatrix}u_{1} & u_{2} & u_{3} \\ u_{1} & u_{1} + u_{2} & v_{3}\end{pmatrix}}
{V=\begin{pmatrix}0 & 0 \\ 0 & 0 \\ 0 & 0\end{pmatrix}}
\subsection{Family \texorpdfstring{$F_{1,4}$}{F1,4} \texorpdfstring{(superdimension $(2\mid4)$)}{(superdimension 2|4)}}
\derivcase{case:1-4-1}{F_{1,4}^{(1)}}
{A=\begin{pmatrix}a_{1} & a_{2} \\ 0 & b_{2}\end{pmatrix}}
{B=\begin{pmatrix}c_{1} & c_{2} & c_{3} & c_{4} \\ 0 & a_{1} + c_{1} & c_{2} & c_{3} \\ 0 & 0 & 2 a_{1} + c_{1} & c_{2} \\ 0 & 0 & 0 & 3 a_{1} + c_{1}\end{pmatrix}}
{U=\begin{pmatrix}u_{1} & u_{2} & u_{3} & u_{4} \\ 0 & 0 & 0 & v_{4}\end{pmatrix}}
{V=\begin{pmatrix}0 & w_{2} \\ 0 & 0 \\ 0 & 0 \\ 0 & 0\end{pmatrix}}
\derivcase{case:1-4-2}{F_{1,4}^{(2)}}
{A=\begin{pmatrix}a_{1} & a_{2} \\ 0 & 2 c_{1}\end{pmatrix}}
{B=\begin{pmatrix}c_{1} & c_{2} & c_{3} & c_{4} \\ 0 & a_{1} + c_{1} & c_{2} & c_{3} \\ 0 & 0 & 2 a_{1} + c_{1} & c_{2} \\ 0 & 0 & 0 & 3 a_{1} + c_{1}\end{pmatrix}}
{U=\begin{pmatrix}u_{1} & u_{2} & u_{3} & u_{4} \\ 0 & 0 & 0 & 0\end{pmatrix}}
{V=\begin{pmatrix}0 & w_{2} \\ 0 & u_{1} \\ 0 & 0 \\ 0 & 0\end{pmatrix}}
\derivcase{case:1-4-3}{F_{1,4}^{(3)}}
{A=\begin{pmatrix}a_{1} & a_{2} \\ 0 & 2 \left(a_{1} + c_{1}\right)\end{pmatrix}}
{B=\begin{pmatrix}c_{1} & c_{2} & 0 & c_{4} \\ 0 & a_{1} + c_{1} & c_{2} & 0 \\ 0 & 0 & 2 a_{1} + c_{1} & c_{2} \\ 0 & 0 & 0 & 3 a_{1} + c_{1}\end{pmatrix}}
{U=\begin{pmatrix}u_{1} & u_{2} & u_{3} & u_{4} \\ 0 & 0 & 0 & 0\end{pmatrix}}
{V=\begin{pmatrix}0 & w_{2} \\ 0 & - u_{3} \\ 0 & u_{2} \\ 0 & - u_{1}\end{pmatrix}}
\derivcase{case:1-4-4}{F_{1,4}^{(4)}}
{A=\begin{pmatrix}a_{1} & a_{2} \\ b_{1} & a_{1} + a_{2} - b_{1}\end{pmatrix}}
{B=\begin{pmatrix}c_{1} & c_{2} & c_{3} & c_{4} \\ 0 & a_{1} + a_{2} + c_{1} & c_{2} & c_{3} \\ 0 & 0 & 2 a_{1} + 2 a_{2} + c_{1} & c_{2} \\ 0 & 0 & 0 & 3 a_{1} + 3 a_{2} + c_{1}\end{pmatrix}}
{U=\begin{pmatrix}u_{1} & u_{2} & u_{3} & u_{4} \\ u_{1} & u_{2} & u_{3} & v_{4}\end{pmatrix}}
{V=\begin{pmatrix}w_{1} & - w_{1} \\ 0 & 0 \\ 0 & 0 \\ 0 & 0\end{pmatrix}}
\derivcase{case:1-4-5}{F_{1,4}^{(5)}}
{A=\begin{pmatrix}a_{1} & 0 \\ 0 & 2 a_{1}\end{pmatrix}}
{B=\begin{pmatrix}c_{1} & c_{2} & c_{3} & c_{4} \\ 0 & a_{1} + c_{1} & c_{2} & c_{3} \\ 0 & 0 & 2 a_{1} + c_{1} & c_{2} \\ 0 & 0 & 0 & 3 a_{1} + c_{1}\end{pmatrix}}
{U=\begin{pmatrix}u_{1} & u_{2} & u_{3} & u_{4} \\ 0 & u_{1} & u_{2} & v_{4}\end{pmatrix}}
{V=\begin{pmatrix}0 & 0 \\ 0 & 0 \\ 0 & 0 \\ 0 & 0\end{pmatrix}}
\derivcase{case:1-4-6}{F_{1,4}^{(6)}}
{A=\begin{pmatrix}a_{1} & a_{2} \\ 0 & 3 a_{1}\end{pmatrix}}
{B=\begin{pmatrix}c_{1} & c_{2} & c_{3} & c_{4} \\ 0 & a_{1} + c_{1} & c_{2} & a_{2} + c_{3} \\ 0 & 0 & 2 a_{1} + c_{1} & c_{2} \\ 0 & 0 & 0 & 3 a_{1} + c_{1}\end{pmatrix}}
{U=\begin{pmatrix}u_{1} & u_{2} & u_{3} & u_{4} \\ 0 & 0 & u_{1} & v_{4}\end{pmatrix}}
{V=\begin{pmatrix}0 & 0 \\ 0 & 0 \\ 0 & 0 \\ 0 & 0\end{pmatrix}}
\derivcase{case:1-4-7}{F_{1,4}^{(7)}}
{A=\begin{pmatrix}a_{1} & 0 \\ - a_{1} & 2 a_{1}\end{pmatrix}}
{B=\begin{pmatrix}c_{1} & c_{2} & c_{3} & c_{4} \\ 0 & a_{1} + c_{1} & c_{2} & c_{3} \\ 0 & 0 & 2 a_{1} + c_{1} & c_{2} \\ 0 & 0 & 0 & 3 a_{1} + c_{1}\end{pmatrix}}
{U=\begin{pmatrix}u_{1} & u_{2} & u_{3} & u_{4} \\ u_{1} & u_{1} + u_{2} & u_{2} + u_{3} & v_{4}\end{pmatrix}}
{V=\begin{pmatrix}0 & 0 \\ 0 & 0 \\ 0 & 0 \\ 0 & 0\end{pmatrix}}
\derivcase{case:1-4-8}{F_{1,4}^{(8)}}
{A=\begin{pmatrix}a_{1} & a_{2} \\ - 2 a_{1} - 3 a_{2} & 3 a_{1} + 4 a_{2}\end{pmatrix}}
{B=\begin{pmatrix}c_{1} & c_{2} & c_{3} & c_{4} \\ 0 & a_{1} + a_{2} + c_{1} & c_{2} & a_{2} + c_{3} \\ 0 & 0 & 2 a_{1} + 2 a_{2} + c_{1} & c_{2} \\ 0 & 0 & 0 & 3 a_{1} + 3 a_{2} + c_{1}\end{pmatrix}}
{U=\begin{pmatrix}u_{1} & u_{2} & u_{3} & u_{4} \\ u_{1} & u_{2} & u_{1} + u_{3} & v_{4}\end{pmatrix}}
{V=\begin{pmatrix}0 & 0 \\ 0 & 0 \\ 0 & 0 \\ 0 & 0\end{pmatrix}}
\derivcase{case:1-4-9}{F_{1,4}^{(9)}}
{A=\begin{pmatrix}a_{1} & - \tfrac{a_{1}}{2} \\ 0 & \tfrac{a_{1}}{2}\end{pmatrix}}
{B=\begin{pmatrix}c_{1} & c_{2} & c_{3} & c_{4} \\ 0 & \tfrac{a_{1} + 2 c_{1}}{2} & - \tfrac{a_{1} - 2 c_{2}}{2} & - a_{1} + c_{3} \\ 0 & 0 & a_{1} + c_{1} & - a_{1} + c_{2} \\ 0 & 0 & 0 & \tfrac{3 a_{1} + 2 c_{1}}{2}\end{pmatrix}}
{U=\begin{pmatrix}u_{1} & u_{2} & u_{3} & u_{4} \\ u_{1} & u_{1} + u_{2} & 2 u_{1} + u_{2} + u_{3} & v_{4}\end{pmatrix}}
{V=\begin{pmatrix}0 & 0 \\ 0 & 0 \\ 0 & 0 \\ 0 & 0\end{pmatrix}}
\subsection{Family \texorpdfstring{$F_{2,2}$}{F2,2} \texorpdfstring{(superdimension $(3\mid2)$)}{(superdimension 3|2)}}
\derivcase{case:2-2-1}{F_{2,2}^{(1)}}
{A=\begin{pmatrix}a_{1} & a_{2} & a_{3} \\ 0 & b_{2} & b_{3} \\ 0 & 0 & a_{1} + b_{2}\end{pmatrix}}
{B=\begin{pmatrix}c_{1} & c_{2} \\ 0 & a_{1} + c_{1}\end{pmatrix}}
{U=\begin{pmatrix}u_{1} & u_{2} \\ v_{1} & v_{2} \\ 0 & v_{1}\end{pmatrix}}
{V=\begin{pmatrix}0 & w_{2} & w_{3} \\ 0 & 0 & w_{2}\end{pmatrix}}
\derivcase{case:2-2-2}{F_{2,2}^{(2)}}
{A=\begin{pmatrix}a_{1} & a_{2} & a_{3} \\ 0 & 2 c_{1} & c_{2} \\ 0 & 0 & a_{1} + 2 c_{1}\end{pmatrix}}
{B=\begin{pmatrix}c_{1} & c_{2} \\ 0 & a_{1} + c_{1}\end{pmatrix}}
{U=\begin{pmatrix}0 & u_{2} \\ 0 & 2 w_{1} \\ 0 & 0\end{pmatrix}}
{V=\begin{pmatrix}w_{1} & w_{2} & w_{3} \\ 0 & 0 & \tfrac{u_{2} + 2 w_{2}}{2}\end{pmatrix}}
\derivcase{case:2-2-3}{F_{2,2}^{(3)}}
{A=\begin{pmatrix}a_{1} & a_{2} & a_{3} \\ 0 & - a_{1} + 2 c_{1} & b_{3} \\ 0 & 0 & 2 c_{1}\end{pmatrix}}
{B=\begin{pmatrix}c_{1} & c_{2} \\ 0 & a_{1} + c_{1}\end{pmatrix}}
{U=\begin{pmatrix}u_{1} & u_{2} \\ 0 & v_{2} \\ 0 & 0\end{pmatrix}}
{V=\begin{pmatrix}0 & w_{2} & w_{3} \\ 0 & 0 & u_{1} + w_{2}\end{pmatrix}}
\derivcase{case:2-2-4}{F_{2,2}^{(4)}}
{A=\begin{pmatrix}a_{1} & a_{2} & a_{3} \\ b_{1} & a_{1} + a_{2} - b_{1} & b_{3} \\ 0 & 0 & 2 a_{1} + a_{2} - b_{1}\end{pmatrix}}
{B=\begin{pmatrix}c_{1} & c_{2} \\ 0 & a_{1} + a_{2} + c_{1}\end{pmatrix}}
{U=\begin{pmatrix}u_{1} & u_{2} \\ v_{1} & v_{2} \\ 0 & - u_{1} + v_{1}\end{pmatrix}}
{V=\begin{pmatrix}w_{1} & - w_{1} & w_{3} \\ 0 & 0 & - w_{1}\end{pmatrix}}
\derivcase{case:2-2-5}{F_{2,2}^{(5)}}
{A=\begin{pmatrix}a_{1} & a_{2} & a_{3} \\ 2 a_{1} + a_{2} - 2 c_{1} & - a_{1} + 2 c_{1} & b_{3} \\ 0 & 0 & 2 c_{1}\end{pmatrix}}
{B=\begin{pmatrix}c_{1} & c_{2} \\ 0 & a_{1} + a_{2} + c_{1}\end{pmatrix}}
{U=\begin{pmatrix}u_{1} & u_{2} \\ u_{1} & v_{2} \\ 0 & 0\end{pmatrix}}
{V=\begin{pmatrix}w_{1} & - w_{1} & w_{3} \\ 0 & 0 & u_{1} - w_{1}\end{pmatrix}}
\subsection{Family \texorpdfstring{$F_{2,3}$}{F2,3} \texorpdfstring{(superdimension $(3\mid3)$)}{(superdimension 3|3)}}
\derivcase{case:2-3-1}{F_{2,3}^{(1)}}
{A=\begin{pmatrix}a_{1} & a_{2} & a_{3} \\ 0 & b_{2} & b_{3} \\ 0 & 0 & a_{1} + b_{2}\end{pmatrix}}
{B=\begin{pmatrix}c_{1} & c_{2} & c_{3} \\ 0 & a_{1} + c_{1} & c_{2} \\ 0 & 0 & 2 a_{1} + c_{1}\end{pmatrix}}
{U=\begin{pmatrix}u_{1} & u_{2} & u_{3} \\ 0 & v_{2} & v_{3} \\ 0 & 0 & v_{2}\end{pmatrix}}
{V=\begin{pmatrix}0 & w_{2} & w_{3} \\ 0 & 0 & w_{2} \\ 0 & 0 & 0\end{pmatrix}}
\derivcase{case:2-3-2}{F_{2,3}^{(2)}}
{A=\begin{pmatrix}a_{1} & a_{2} & a_{3} \\ b_{1} & a_{1} + a_{2} - b_{1} & b_{3} \\ 0 & 0 & 2 a_{1} + a_{2} - b_{1}\end{pmatrix}}
{B=\begin{pmatrix}c_{1} & c_{2} & c_{3} \\ 0 & a_{1} + a_{2} + c_{1} & c_{2} \\ 0 & 0 & 2 a_{1} + 2 a_{2} + c_{1}\end{pmatrix}}
{U=\begin{pmatrix}u_{1} & u_{2} & u_{3} \\ u_{1} & v_{2} & v_{3} \\ 0 & 0 & - u_{2} + v_{2}\end{pmatrix}}
{V=\begin{pmatrix}w_{1} & - w_{1} & w_{3} \\ 0 & 0 & - w_{1} \\ 0 & 0 & 0\end{pmatrix}}
\derivcase{case:2-3-3}{F_{2,3}^{(3)}}
{A=\begin{pmatrix}a_{1} & a_{2} & a_{3} \\ 0 & 2 a_{1} & b_{3} \\ 0 & 0 & 3 a_{1}\end{pmatrix}}
{B=\begin{pmatrix}c_{1} & c_{2} & c_{3} \\ 0 & a_{1} + c_{1} & a_{2} + c_{2} \\ 0 & 0 & 2 a_{1} + c_{1}\end{pmatrix}}
{U=\begin{pmatrix}u_{1} & u_{2} & u_{3} \\ 0 & v_{2} & v_{3} \\ 0 & 0 & - u_{1} + v_{2}\end{pmatrix}}
{V=\begin{pmatrix}0 & 0 & w_{3} \\ 0 & 0 & 0 \\ 0 & 0 & 0\end{pmatrix}}
\derivcase{case:2-3-4}{F_{2,3}^{(4)}}
{A=\begin{pmatrix}a_{1} & a_{2} & a_{3} \\ 0 & a_{1} - a_{2} & c_{2} \\ 0 & 0 & 2 a_{1} - a_{2}\end{pmatrix}}
{B=\begin{pmatrix}c_{1} & c_{2} & c_{3} \\ 0 & a_{1} + c_{1} & a_{3} + c_{2} \\ 0 & 0 & 2 a_{1} - a_{2} + c_{1}\end{pmatrix}}
{U=\begin{pmatrix}0 & u_{2} & u_{3} \\ 0 & v_{2} & v_{3} \\ 0 & 0 & u_{2} + v_{2}\end{pmatrix}}
{V=\begin{pmatrix}0 & w_{2} & 0 \\ 0 & 0 & w_{2} \\ 0 & 0 & 0\end{pmatrix}}
\derivcase{case:2-3-5}{F_{2,3}^{(5)}}
{A=\begin{pmatrix}a_{1} & a_{2} & a_{3} \\ - a_{1} - 2 a_{2} & 2 a_{1} + 3 a_{2} & b_{3} \\ 0 & 0 & 3 \left(a_{1} + a_{2}\right)\end{pmatrix}}
{B=\begin{pmatrix}c_{1} & c_{2} & c_{3} \\ 0 & a_{1} + a_{2} + c_{1} & a_{2} + c_{2} \\ 0 & 0 & 2 a_{1} + 2 a_{2} + c_{1}\end{pmatrix}}
{U=\begin{pmatrix}u_{1} & u_{2} & u_{3} \\ u_{1} & v_{2} & v_{3} \\ 0 & 0 & - u_{1} - u_{2} + v_{2}\end{pmatrix}}
{V=\begin{pmatrix}0 & 0 & w_{3} \\ 0 & 0 & 0 \\ 0 & 0 & 0\end{pmatrix}}
\derivcase{case:2-3-6}{F_{2,3}^{(6)}}
{A=\begin{pmatrix}a_{1} & a_{2} & a_{3} \\ 0 & a_{1} + a_{2} & a_{3} + c_{2} \\ 0 & 0 & 2 a_{1} + a_{2}\end{pmatrix}}
{B=\begin{pmatrix}c_{1} & c_{2} & c_{3} \\ 0 & a_{1} + a_{2} + c_{1} & a_{3} + c_{2} \\ 0 & 0 & 2 a_{1} + a_{2} + c_{1}\end{pmatrix}}
{U=\begin{pmatrix}0 & u_{2} & u_{3} \\ 0 & v_{2} & v_{3} \\ 0 & 0 & v_{2}\end{pmatrix}}
{V=\begin{pmatrix}w_{1} & - w_{1} & 0 \\ 0 & 0 & - w_{1} \\ 0 & 0 & 0\end{pmatrix}}
\derivcase{case:2-3-7}{F_{2,3}^{(7)}}
{A=\begin{pmatrix}a_{1} & a_{2} & a_{3} \\ - 2 a_{2} & a_{1} + 3 a_{2} & 2 a_{3} + c_{2} \\ 0 & 0 & 2 a_{1} + 3 a_{2}\end{pmatrix}}
{B=\begin{pmatrix}c_{1} & c_{2} & c_{3} \\ 0 & a_{1} + 2 a_{2} + c_{1} & a_{3} + c_{2} \\ 0 & 0 & 2 a_{1} + 3 a_{2} + c_{1}\end{pmatrix}}
{U=\begin{pmatrix}0 & u_{2} & u_{3} \\ 0 & v_{2} & v_{3} \\ 0 & 0 & - u_{2} + v_{2}\end{pmatrix}}
{V=\begin{pmatrix}w_{1} & - \tfrac{w_{1}}{2} & 0 \\ 0 & 0 & - \tfrac{w_{1}}{2} \\ 0 & 0 & 0\end{pmatrix}}
\derivcase{case:2-3-8}{F_{2,3}^{(8)}}
{A=\begin{pmatrix}a_{1} & a_{2} & a_{3} \\ 0 & - a_{1} + 2 c_{1} & b_{3} \\ 0 & 0 & 2 c_{1}\end{pmatrix}}
{B=\begin{pmatrix}c_{1} & c_{2} & c_{3} \\ 0 & a_{1} + c_{1} & c_{2} \\ 0 & 0 & 2 a_{1} + c_{1}\end{pmatrix}}
{U=\begin{pmatrix}u_{1} & u_{2} & u_{3} \\ 0 & 0 & v_{3} \\ 0 & 0 & 0\end{pmatrix}}
{V=\begin{pmatrix}0 & w_{2} & w_{3} \\ 0 & 0 & u_{1} + w_{2} \\ 0 & 0 & 0\end{pmatrix}}
\derivcase{case:2-3-9}{F_{2,3}^{(9)}}
{A=\begin{pmatrix}a_{1} & a_{2} & a_{3} \\ 2 a_{1} + a_{2} - 2 c_{1} & - a_{1} + 2 c_{1} & b_{3} \\ 0 & 0 & 2 c_{1}\end{pmatrix}}
{B=\begin{pmatrix}c_{1} & c_{2} & c_{3} \\ 0 & a_{1} + a_{2} + c_{1} & c_{2} \\ 0 & 0 & 2 a_{1} + 2 a_{2} + c_{1}\end{pmatrix}}
{U=\begin{pmatrix}u_{1} & u_{2} & u_{3} \\ u_{1} & u_{2} & v_{3} \\ 0 & 0 & 0\end{pmatrix}}
{V=\begin{pmatrix}w_{1} & - w_{1} & w_{3} \\ 0 & 0 & u_{1} - w_{1} \\ 0 & 0 & 0\end{pmatrix}}
\derivcase{case:2-3-10}{F_{2,3}^{(10)}}
{A=\begin{pmatrix}a_{1} & a_{2} & a_{3} \\ 0 & 2 a_{1} & b_{3} \\ 0 & 0 & 3 a_{1}\end{pmatrix}}
{B=\begin{pmatrix}\tfrac{3 a_{1}}{2} & c_{2} & c_{3} \\ 0 & \tfrac{5 a_{1}}{2} & a_{2} + c_{2} \\ 0 & 0 & \tfrac{7 a_{1}}{2}\end{pmatrix}}
{U=\begin{pmatrix}u_{1} & u_{2} & u_{3} \\ 0 & u_{1} + 2 w_{2} & v_{3} \\ 0 & 0 & 2 w_{2}\end{pmatrix}}
{V=\begin{pmatrix}0 & w_{2} & w_{3} \\ 0 & 0 & u_{1} + w_{2} \\ 0 & 0 & 0\end{pmatrix}}
\derivcase{case:2-3-11}{F_{2,3}^{(11)}}
{A=\begin{pmatrix}a_{1} & a_{2} & a_{3} \\ - a_{1} - 2 a_{2} & 2 a_{1} + 3 a_{2} & b_{3} \\ 0 & 0 & 3 \left(a_{1} + a_{2}\right)\end{pmatrix}}
{B=\begin{pmatrix}\tfrac{3 \left(a_{1} + a_{2}\right)}{2} & c_{2} & c_{3} \\ 0 & \tfrac{5 \left(a_{1} + a_{2}\right)}{2} & a_{2} + c_{2} \\ 0 & 0 & \tfrac{7 \left(a_{1} + a_{2}\right)}{2}\end{pmatrix}}
{U=\begin{pmatrix}u_{1} & u_{2} & u_{3} \\ u_{1} & u_{1} + u_{2} - 2 w_{1} & v_{3} \\ 0 & 0 & - 2 w_{1}\end{pmatrix}}
{V=\begin{pmatrix}w_{1} & - w_{1} & w_{3} \\ 0 & 0 & u_{1} - w_{1} \\ 0 & 0 & 0\end{pmatrix}}
\derivcase{case:2-3-12}{F_{2,3}^{(12)}}
{A=\begin{pmatrix}a_{1} & a_{2} & a_{3} \\ 0 & a_{1} + 2 c_{1} & b_{3} \\ 0 & 0 & 2 \left(a_{1} + c_{1}\right)\end{pmatrix}}
{B=\begin{pmatrix}c_{1} & c_{2} & 0 \\ 0 & a_{1} + c_{1} & c_{2} \\ 0 & 0 & 2 a_{1} + c_{1}\end{pmatrix}}
{U=\begin{pmatrix}0 & u_{2} & u_{3} \\ 0 & 0 & 0 \\ 0 & 0 & 0\end{pmatrix}}
{V=\begin{pmatrix}0 & w_{2} & w_{3} \\ 0 & 0 & - u_{3} + w_{2} \\ 0 & 0 & u_{2}\end{pmatrix}}
\derivcase{case:2-3-13}{F_{2,3}^{(13)}}
{A=\begin{pmatrix}a_{1} & a_{2} & a_{3} \\ - a_{2} - 2 c_{1} & a_{1} + 2 a_{2} + 2 c_{1} & b_{3} \\ 0 & 0 & 2 \left(a_{1} + a_{2} + c_{1}\right)\end{pmatrix}}
{B=\begin{pmatrix}c_{1} & c_{2} & 0 \\ 0 & a_{1} + a_{2} + c_{1} & c_{2} \\ 0 & 0 & 2 a_{1} + 2 a_{2} + c_{1}\end{pmatrix}}
{U=\begin{pmatrix}0 & u_{2} & u_{3} \\ 0 & u_{2} & u_{3} \\ 0 & 0 & 0\end{pmatrix}}
{V=\begin{pmatrix}w_{1} & - w_{1} & w_{3} \\ 0 & 0 & - u_{3} - w_{1} \\ 0 & 0 & u_{2}\end{pmatrix}}
\derivcase{case:2-3-14}{F_{2,3}^{(14)}}
{A=\begin{pmatrix}a_{1} & a_{2} & a_{3} \\ 0 & 2 c_{1} & c_{2} \\ 0 & 0 & a_{1} + 2 c_{1}\end{pmatrix}}
{B=\begin{pmatrix}c_{1} & c_{2} & c_{3} \\ 0 & a_{1} + c_{1} & c_{2} \\ 0 & 0 & 2 a_{1} + c_{1}\end{pmatrix}}
{U=\begin{pmatrix}u_{1} & u_{2} & u_{3} \\ 0 & 2 w_{1} & 0 \\ 0 & 0 & 2 w_{1}\end{pmatrix}}
{V=\begin{pmatrix}w_{1} & w_{2} & w_{3} \\ 0 & u_{1} & \tfrac{u_{2} + 2 w_{2}}{2} \\ 0 & 0 & \tfrac{3 u_{1}}{2}\end{pmatrix}}
\derivcase{case:2-3-15}{F_{2,3}^{(15)}}
{A=\begin{pmatrix}a_{1} & a_{2} & a_{3} \\ 0 & a_{1} - a_{2} & c_{2} \\ 0 & 0 & 2 a_{1} - a_{2}\end{pmatrix}}
{B=\begin{pmatrix}\tfrac{a_{1} - a_{2}}{2} & c_{2} & c_{3} \\ 0 & \tfrac{3 a_{1} - a_{2}}{2} & a_{3} + c_{2} \\ 0 & 0 & \tfrac{5 a_{1} - 3 a_{2}}{2}\end{pmatrix}}
{U=\begin{pmatrix}0 & u_{2} & u_{3} \\ 0 & 2 w_{1} & 2 w_{3} \\ 0 & 0 & u_{2} + 2 w_{1}\end{pmatrix}}
{V=\begin{pmatrix}w_{1} & w_{2} & w_{3} \\ 0 & 0 & \tfrac{u_{2} + 2 w_{2}}{2} \\ 0 & 0 & 0\end{pmatrix}}
\derivcase{case:2-3-16}{F_{2,3}^{(16)}}
{A=\begin{pmatrix}0 & a_{2} & a_{3} \\ 0 & 2 c_{1} & c_{2} - 2 c_{3} \\ 0 & 0 & 2 c_{1}\end{pmatrix}}
{B=\begin{pmatrix}c_{1} & c_{2} & c_{3} \\ 0 & c_{1} & c_{2} \\ 0 & 0 & c_{1}\end{pmatrix}}
{U=\begin{pmatrix}0 & u_{2} & u_{3} \\ 0 & 0 & 0 \\ 0 & 0 & 0\end{pmatrix}}
{V=\begin{pmatrix}0 & w_{2} & w_{3} \\ 0 & 0 & \tfrac{u_{2} - 2 u_{3} + 2 w_{2}}{2} \\ 0 & 0 & u_{2}\end{pmatrix}}
\subsection{Family \texorpdfstring{$F_{3,2}$}{F3,2} \texorpdfstring{(superdimension $(4\mid2)$)}{(superdimension 4|2)}}
\derivcase{case:3-2-1}{F_{3,2}^{(1)}}
{A=\begin{pmatrix}a_{1} & a_{2} & a_{3} & a_{4} \\ 0 & b_{2} & b_{3} & b_{4} \\ 0 & 0 & a_{1} + b_{2} & b_{3} \\ 0 & 0 & 0 & 2 a_{1} + b_{2}\end{pmatrix}}
{B=\begin{pmatrix}c_{1} & c_{2} \\ 0 & a_{1} + c_{1}\end{pmatrix}}
{U=\begin{pmatrix}u_{1} & u_{2} \\ v_{1} & v_{2} \\ 0 & v_{1} \\ 0 & 0\end{pmatrix}}
{V=\begin{pmatrix}0 & 0 & w_{3} & w_{4} \\ 0 & 0 & 0 & w_{3}\end{pmatrix}}
\derivcase{case:3-2-2}{F_{3,2}^{(2)}}
{A=\begin{pmatrix}a_{1} & a_{2} & a_{3} & a_{4} \\ 0 & 2 c_{1} & c_{2} & 0 \\ 0 & 0 & a_{1} + 2 c_{1} & c_{2} \\ 0 & 0 & 0 & 2 \left(a_{1} + c_{1}\right)\end{pmatrix}}
{B=\begin{pmatrix}c_{1} & c_{2} \\ 0 & a_{1} + c_{1}\end{pmatrix}}
{U=\begin{pmatrix}0 & u_{2} \\ 0 & 0 \\ 0 & 0 \\ 0 & 0\end{pmatrix}}
{V=\begin{pmatrix}0 & - u_{2} & w_{3} & w_{4} \\ 0 & 0 & - \tfrac{u_{2}}{2} & w_{3}\end{pmatrix}}
\derivcase{case:3-2-3}{F_{3,2}^{(3)}}
{A=\begin{pmatrix}a_{1} & a_{2} & a_{3} & a_{4} \\ 0 & - a_{1} + 2 c_{1} & c_{2} & b_{4} \\ 0 & 0 & 2 c_{1} & c_{2} \\ 0 & 0 & 0 & a_{1} + 2 c_{1}\end{pmatrix}}
{B=\begin{pmatrix}c_{1} & c_{2} \\ 0 & a_{1} + c_{1}\end{pmatrix}}
{U=\begin{pmatrix}u_{1} & u_{2} \\ 0 & 0 \\ 0 & 0 \\ 0 & 0\end{pmatrix}}
{V=\begin{pmatrix}0 & - \tfrac{3 u_{1}}{2} & w_{3} & w_{4} \\ 0 & 0 & - \tfrac{u_{1}}{2} & \tfrac{u_{2} + 2 w_{3}}{2}\end{pmatrix}}
\derivcase{case:3-2-4}{F_{3,2}^{(4)}}
{A=\begin{pmatrix}a_{1} & a_{2} & a_{3} & a_{4} \\ 0 & - 2 \left(a_{1} - c_{1}\right) & b_{3} & b_{4} \\ 0 & 0 & - a_{1} + 2 c_{1} & b_{3} \\ 0 & 0 & 0 & 2 c_{1}\end{pmatrix}}
{B=\begin{pmatrix}c_{1} & c_{2} \\ 0 & a_{1} + c_{1}\end{pmatrix}}
{U=\begin{pmatrix}u_{1} & u_{2} \\ 0 & v_{2} \\ 0 & 0 \\ 0 & 0\end{pmatrix}}
{V=\begin{pmatrix}0 & 0 & w_{3} & w_{4} \\ 0 & 0 & 0 & u_{1} + w_{3}\end{pmatrix}}
\derivcase{case:3-2-5}{F_{3,2}^{(5)}}
{A=\begin{pmatrix}a_{1} & a_{2} & a_{3} & a_{4} \\ 0 & a_{1} + a_{2} & b_{3} & b_{4} \\ 0 & 0 & 2 a_{1} + a_{2} & b_{3} \\ 0 & 0 & 0 & 3 a_{1} + a_{2}\end{pmatrix}}
{B=\begin{pmatrix}c_{1} & c_{2} \\ 0 & a_{1} + a_{2} + c_{1}\end{pmatrix}}
{U=\begin{pmatrix}u_{1} & u_{2} \\ v_{1} & v_{2} \\ 0 & - u_{1} + v_{1} \\ 0 & 0\end{pmatrix}}
{V=\begin{pmatrix}0 & 0 & 0 & w_{4} \\ 0 & 0 & 0 & 0\end{pmatrix}}
\derivcase{case:3-2-6}{F_{3,2}^{(6)}}
{A=\begin{pmatrix}a_{1} & a_{2} & a_{3} & a_{4} \\ 0 & a_{1} + a_{2} & b_{3} & b_{4} \\ 0 & 0 & 2 a_{1} + a_{2} & b_{3} \\ 0 & 0 & 0 & 3 a_{1} + a_{2}\end{pmatrix}}
{B=\begin{pmatrix}\tfrac{3 a_{1} + a_{2}}{2} & c_{2} \\ 0 & \tfrac{5 a_{1} + 3 a_{2}}{2}\end{pmatrix}}
{U=\begin{pmatrix}u_{1} & u_{2} \\ u_{1} + w_{3} & v_{2} \\ 0 & w_{3} \\ 0 & 0\end{pmatrix}}
{V=\begin{pmatrix}0 & 0 & w_{3} & w_{4} \\ 0 & 0 & 0 & u_{1} + w_{3}\end{pmatrix}}
\subsection{Family \texorpdfstring{$F_{4,2}$}{F4,2} \texorpdfstring{(superdimension $(5\mid2)$)}{(superdimension 5|2)}}
\derivcase{case:4-2-1}{F_{4,2}^{(1)}}
{A=\begin{pmatrix}a_{1} & a_{2} & a_{3} & a_{4} & a_{5} \\ 0 & b_{2} & b_{3} & b_{4} & b_{5} \\ 0 & 0 & a_{1} + b_{2} & b_{3} & b_{4} \\ 0 & 0 & 0 & 2 a_{1} + b_{2} & b_{3} \\ 0 & 0 & 0 & 0 & 3 a_{1} + b_{2}\end{pmatrix}}
{B=\begin{pmatrix}c_{1} & c_{2} \\ 0 & a_{1} + c_{1}\end{pmatrix}}
{U=\begin{pmatrix}u_{1} & u_{2} \\ v_{1} & v_{2} \\ 0 & v_{1} \\ 0 & 0 \\ 0 & 0\end{pmatrix}}
{V=\begin{pmatrix}0 & 0 & 0 & w_{4} & w_{5} \\ 0 & 0 & 0 & 0 & w_{4}\end{pmatrix}}
\derivcase{case:4-2-2}{F_{4,2}^{(2)}}
{A=\begin{pmatrix}a_{1} & a_{2} & a_{3} & a_{4} & a_{5} \\ 0 & - a_{1} + 2 c_{1} & c_{2} & 0 & b_{5} \\ 0 & 0 & 2 c_{1} & c_{2} & 0 \\ 0 & 0 & 0 & a_{1} + 2 c_{1} & c_{2} \\ 0 & 0 & 0 & 0 & 2 \left(a_{1} + c_{1}\right)\end{pmatrix}}
{B=\begin{pmatrix}c_{1} & c_{2} \\ 0 & a_{1} + c_{1}\end{pmatrix}}
{U=\begin{pmatrix}u_{1} & u_{2} \\ 0 & 0 \\ 0 & 0 \\ 0 & 0 \\ 0 & 0\end{pmatrix}}
{V=\begin{pmatrix}0 & - \tfrac{3 u_{1}}{2} & - u_{2} & w_{4} & w_{5} \\ 0 & 0 & - \tfrac{u_{1}}{2} & - \tfrac{u_{2}}{2} & w_{4}\end{pmatrix}}
\derivcase{case:4-2-3}{F_{4,2}^{(3)}}
{A=\begin{pmatrix}a_{1} & a_{2} & a_{3} & a_{4} & a_{5} \\ 0 & - 2 \left(a_{1} - c_{1}\right) & c_{2} & b_{4} & b_{5} \\ 0 & 0 & - a_{1} + 2 c_{1} & c_{2} & b_{4} \\ 0 & 0 & 0 & 2 c_{1} & c_{2} \\ 0 & 0 & 0 & 0 & a_{1} + 2 c_{1}\end{pmatrix}}
{B=\begin{pmatrix}c_{1} & c_{2} \\ 0 & a_{1} + c_{1}\end{pmatrix}}
{U=\begin{pmatrix}u_{1} & u_{2} \\ 0 & 0 \\ 0 & 0 \\ 0 & 0 \\ 0 & 0\end{pmatrix}}
{V=\begin{pmatrix}0 & 0 & - \tfrac{3 u_{1}}{2} & w_{4} & w_{5} \\ 0 & 0 & 0 & - \tfrac{u_{1}}{2} & \tfrac{u_{2} + 2 w_{4}}{2}\end{pmatrix}}
\derivcase{case:4-2-4}{F_{4,2}^{(4)}}
{A=\begin{pmatrix}a_{1} & a_{2} & a_{3} & a_{4} & a_{5} \\ 0 & - 3 a_{1} + 2 c_{1} & b_{3} & b_{4} & b_{5} \\ 0 & 0 & - 2 \left(a_{1} - c_{1}\right) & b_{3} & b_{4} \\ 0 & 0 & 0 & - a_{1} + 2 c_{1} & b_{3} \\ 0 & 0 & 0 & 0 & 2 c_{1}\end{pmatrix}}
{B=\begin{pmatrix}c_{1} & c_{2} \\ 0 & a_{1} + c_{1}\end{pmatrix}}
{U=\begin{pmatrix}u_{1} & u_{2} \\ 0 & v_{2} \\ 0 & 0 \\ 0 & 0 \\ 0 & 0\end{pmatrix}}
{V=\begin{pmatrix}0 & 0 & 0 & w_{4} & w_{5} \\ 0 & 0 & 0 & 0 & u_{1} + w_{4}\end{pmatrix}}
\derivcase{case:4-2-5}{F_{4,2}^{(5)}}
{A=\begin{pmatrix}a_{1} & a_{2} & a_{3} & a_{4} & a_{5} \\ 0 & a_{1} + a_{2} & b_{3} & b_{4} & b_{5} \\ 0 & 0 & 2 a_{1} + a_{2} & b_{3} & b_{4} \\ 0 & 0 & 0 & 3 a_{1} + a_{2} & b_{3} \\ 0 & 0 & 0 & 0 & 4 a_{1} + a_{2}\end{pmatrix}}
{B=\begin{pmatrix}c_{1} & c_{2} \\ 0 & a_{1} + a_{2} + c_{1}\end{pmatrix}}
{U=\begin{pmatrix}u_{1} & u_{2} \\ v_{1} & v_{2} \\ 0 & - u_{1} + v_{1} \\ 0 & 0 \\ 0 & 0\end{pmatrix}}
{V=\begin{pmatrix}0 & 0 & 0 & 0 & w_{5} \\ 0 & 0 & 0 & 0 & 0\end{pmatrix}}
\derivcase{case:4-2-6}{F_{4,2}^{(6)}}
{A=\begin{pmatrix}a_{1} & a_{2} & a_{3} & a_{4} & a_{5} \\ 0 & a_{1} + a_{2} & b_{3} & b_{4} & b_{5} \\ 0 & 0 & 2 a_{1} + a_{2} & b_{3} & b_{4} \\ 0 & 0 & 0 & 3 a_{1} + a_{2} & b_{3} \\ 0 & 0 & 0 & 0 & 4 a_{1} + a_{2}\end{pmatrix}}
{B=\begin{pmatrix}\tfrac{4 a_{1} + a_{2}}{2} & c_{2} \\ 0 & \tfrac{3 \left(2 a_{1} + a_{2}\right)}{2}\end{pmatrix}}
{U=\begin{pmatrix}u_{1} & u_{2} \\ u_{1} + w_{4} & v_{2} \\ 0 & w_{4} \\ 0 & 0 \\ 0 & 0\end{pmatrix}}
{V=\begin{pmatrix}0 & 0 & 0 & w_{4} & w_{5} \\ 0 & 0 & 0 & 0 & u_{1} + w_{4}\end{pmatrix}}
\derivcase{case:4-2-7}{F_{4,2}^{(7)}}
{A=\begin{pmatrix}a_{1} & a_{2} & a_{3} & a_{4} & a_{5} \\ 0 & 2 a_{1} & b_{3} & b_{4} & b_{5} \\ 0 & 0 & 3 a_{1} & b_{3} & - a_{3} + b_{4} \\ 0 & 0 & 0 & 4 a_{1} & a_{2} + b_{3} \\ 0 & 0 & 0 & 0 & 5 a_{1}\end{pmatrix}}
{B=\begin{pmatrix}c_{1} & c_{2} \\ 0 & a_{1} + c_{1}\end{pmatrix}}
{U=\begin{pmatrix}u_{1} & u_{2} \\ v_{1} & v_{2} \\ 0 & v_{1} \\ 0 & 0 \\ 0 & 0\end{pmatrix}}
{V=\begin{pmatrix}0 & 0 & 0 & w_{4} & w_{5} \\ 0 & 0 & 0 & 0 & w_{4}\end{pmatrix}}
\derivcase{case:4-2-8}{F_{4,2}^{(8)}}
{A=\begin{pmatrix}a_{1} & a_{2} & a_{3} & a_{4} & a_{5} \\ 0 & 2 a_{1} & - a_{2} + c_{2} & b_{4} & b_{5} \\ 0 & 0 & 3 a_{1} & - a_{2} + c_{2} & - a_{3} + b_{4} \\ 0 & 0 & 0 & 4 a_{1} & c_{2} \\ 0 & 0 & 0 & 0 & 5 a_{1}\end{pmatrix}}
{B=\begin{pmatrix}2 a_{1} & c_{2} \\ 0 & 3 a_{1}\end{pmatrix}}
{U=\begin{pmatrix}u_{1} & u_{2} \\ 0 & 3 u_{1} \\ 0 & 0 \\ 0 & 0 \\ 0 & 0\end{pmatrix}}
{V=\begin{pmatrix}0 & 0 & - \tfrac{3 u_{1}}{2} & w_{4} & w_{5} \\ 0 & 0 & 0 & - \tfrac{u_{1}}{2} & \tfrac{u_{2} + 2 w_{4}}{2}\end{pmatrix}}
\derivcase{case:4-2-9}{F_{4,2}^{(9)}}
{A=\begin{pmatrix}a_{1} & a_{2} & a_{3} & a_{4} & a_{5} \\ 0 & 2 a_{1} & b_{3} & b_{4} & b_{5} \\ 0 & 0 & 3 a_{1} & b_{3} & - a_{3} + b_{4} \\ 0 & 0 & 0 & 4 a_{1} & a_{2} + b_{3} \\ 0 & 0 & 0 & 0 & 5 a_{1}\end{pmatrix}}
{B=\begin{pmatrix}\tfrac{5 a_{1}}{2} & c_{2} \\ 0 & \tfrac{7 a_{1}}{2}\end{pmatrix}}
{U=\begin{pmatrix}u_{1} & u_{2} \\ 0 & v_{2} \\ 0 & 0 \\ 0 & 0 \\ 0 & 0\end{pmatrix}}
{V=\begin{pmatrix}0 & 0 & 0 & w_{4} & w_{5} \\ 0 & 0 & 0 & 0 & u_{1} + w_{4}\end{pmatrix}}
\derivcase{case:4-2-10}{F_{4,2}^{(10)}}
{A=\begin{pmatrix}a_{1} & a_{1} & a_{3} & a_{4} & a_{5} \\ 0 & 2 a_{1} & b_{3} & b_{4} & b_{5} \\ 0 & 0 & 3 a_{1} & b_{3} & - a_{3} + b_{4} \\ 0 & 0 & 0 & 4 a_{1} & a_{1} + b_{3} \\ 0 & 0 & 0 & 0 & 5 a_{1}\end{pmatrix}}
{B=\begin{pmatrix}c_{1} & c_{2} \\ 0 & 2 a_{1} + c_{1}\end{pmatrix}}
{U=\begin{pmatrix}u_{1} & u_{2} \\ v_{1} & v_{2} \\ 0 & - u_{1} + v_{1} \\ 0 & 0 \\ 0 & 0\end{pmatrix}}
{V=\begin{pmatrix}0 & 0 & 0 & 0 & w_{5} \\ 0 & 0 & 0 & 0 & 0\end{pmatrix}}
\derivcase{case:4-2-11}{F_{4,2}^{(11)}}
{A=\begin{pmatrix}a_{1} & a_{1} & a_{3} & a_{4} & a_{5} \\ 0 & 2 a_{1} & b_{3} & b_{4} & b_{5} \\ 0 & 0 & 3 a_{1} & b_{3} & - a_{3} + b_{4} \\ 0 & 0 & 0 & 4 a_{1} & a_{1} + b_{3} \\ 0 & 0 & 0 & 0 & 5 a_{1}\end{pmatrix}}
{B=\begin{pmatrix}\tfrac{5 a_{1}}{2} & c_{2} \\ 0 & \tfrac{9 a_{1}}{2}\end{pmatrix}}
{U=\begin{pmatrix}u_{1} & u_{2} \\ u_{1} + w_{4} & v_{2} \\ 0 & w_{4} \\ 0 & 0 \\ 0 & 0\end{pmatrix}}
{V=\begin{pmatrix}0 & 0 & 0 & w_{4} & w_{5} \\ 0 & 0 & 0 & 0 & u_{1} + w_{4}\end{pmatrix}}
\subsection{Family \texorpdfstring{$F_{5,2}$}{F5,2} \texorpdfstring{(superdimension $(6\mid2)$)}{(superdimension 6|2)}}
\derivcase{case:5-2-1}{F_{5,2}^{(1)}}
{A=\begin{pmatrix}a_{1} & a_{2} & a_{3} & a_{4} & a_{5} & a_{6} \\ 0 & b_{2} & b_{3} & b_{4} & b_{5} & b_{6} \\ 0 & 0 & a_{1} + b_{2} & b_{3} & b_{4} & b_{5} \\ 0 & 0 & 0 & 2 a_{1} + b_{2} & b_{3} & b_{4} \\ 0 & 0 & 0 & 0 & 3 a_{1} + b_{2} & b_{3} \\ 0 & 0 & 0 & 0 & 0 & 4 a_{1} + b_{2}\end{pmatrix}}
{B=\begin{pmatrix}c_{1} & c_{2} \\ 0 & a_{1} + c_{1}\end{pmatrix}}
{U=\begin{pmatrix}u_{1} & u_{2} \\ v_{1} & v_{2} \\ 0 & v_{1} \\ 0 & 0 \\ 0 & 0 \\ 0 & 0\end{pmatrix}}
{V=\begin{pmatrix}0 & 0 & 0 & 0 & w_{5} & w_{6} \\ 0 & 0 & 0 & 0 & 0 & w_{5}\end{pmatrix}}
\derivcase{case:5-2-2}{F_{5,2}^{(2)}}
{A=\begin{pmatrix}a_{1} & a_{2} & a_{3} & a_{4} & a_{5} & a_{6} \\ 0 & - 2 \left(a_{1} - c_{1}\right) & c_{2} & 0 & b_{5} & b_{6} \\ 0 & 0 & - a_{1} + 2 c_{1} & c_{2} & 0 & b_{5} \\ 0 & 0 & 0 & 2 c_{1} & c_{2} & 0 \\ 0 & 0 & 0 & 0 & a_{1} + 2 c_{1} & c_{2} \\ 0 & 0 & 0 & 0 & 0 & 2 \left(a_{1} + c_{1}\right)\end{pmatrix}}
{B=\begin{pmatrix}c_{1} & c_{2} \\ 0 & a_{1} + c_{1}\end{pmatrix}}
{U=\begin{pmatrix}u_{1} & u_{2} \\ 0 & 0 \\ 0 & 0 \\ 0 & 0 \\ 0 & 0 \\ 0 & 0\end{pmatrix}}
{V=\begin{pmatrix}0 & 0 & - \tfrac{3 u_{1}}{2} & - u_{2} & w_{5} & w_{6} \\ 0 & 0 & 0 & - \tfrac{u_{1}}{2} & - \tfrac{u_{2}}{2} & w_{5}\end{pmatrix}}
\derivcase{case:5-2-3}{F_{5,2}^{(3)}}
{A=\begin{pmatrix}a_{1} & a_{2} & a_{3} & a_{4} & a_{5} & a_{6} \\ 0 & - 3 a_{1} + 2 c_{1} & c_{2} & b_{4} & b_{5} & b_{6} \\ 0 & 0 & - 2 \left(a_{1} - c_{1}\right) & c_{2} & b_{4} & b_{5} \\ 0 & 0 & 0 & - a_{1} + 2 c_{1} & c_{2} & b_{4} \\ 0 & 0 & 0 & 0 & 2 c_{1} & c_{2} \\ 0 & 0 & 0 & 0 & 0 & a_{1} + 2 c_{1}\end{pmatrix}}
{B=\begin{pmatrix}c_{1} & c_{2} \\ 0 & a_{1} + c_{1}\end{pmatrix}}
{U=\begin{pmatrix}u_{1} & u_{2} \\ 0 & 0 \\ 0 & 0 \\ 0 & 0 \\ 0 & 0 \\ 0 & 0\end{pmatrix}}
{V=\begin{pmatrix}0 & 0 & 0 & - \tfrac{3 u_{1}}{2} & w_{5} & w_{6} \\ 0 & 0 & 0 & 0 & - \tfrac{u_{1}}{2} & \tfrac{u_{2} + 2 w_{5}}{2}\end{pmatrix}}
\derivcase{case:5-2-4}{F_{5,2}^{(4)}}
{A=\begin{pmatrix}a_{1} & a_{2} & a_{3} & a_{4} & a_{5} & a_{6} \\ 0 & - 2 \left(2 a_{1} - c_{1}\right) & b_{3} & b_{4} & b_{5} & b_{6} \\ 0 & 0 & - 3 a_{1} + 2 c_{1} & b_{3} & b_{4} & b_{5} \\ 0 & 0 & 0 & - 2 \left(a_{1} - c_{1}\right) & b_{3} & b_{4} \\ 0 & 0 & 0 & 0 & - a_{1} + 2 c_{1} & b_{3} \\ 0 & 0 & 0 & 0 & 0 & 2 c_{1}\end{pmatrix}}
{B=\begin{pmatrix}c_{1} & c_{2} \\ 0 & a_{1} + c_{1}\end{pmatrix}}
{U=\begin{pmatrix}u_{1} & u_{2} \\ 0 & v_{2} \\ 0 & 0 \\ 0 & 0 \\ 0 & 0 \\ 0 & 0\end{pmatrix}}
{V=\begin{pmatrix}0 & 0 & 0 & 0 & w_{5} & w_{6} \\ 0 & 0 & 0 & 0 & 0 & u_{1} + w_{5}\end{pmatrix}}
\derivcase{case:5-2-5}{F_{5,2}^{(5)}}
{A=\begin{pmatrix}a_{1} & a_{2} & a_{3} & a_{4} & a_{5} & a_{6} \\ 0 & a_{1} + a_{2} & b_{3} & b_{4} & b_{5} & b_{6} \\ 0 & 0 & 2 a_{1} + a_{2} & b_{3} & b_{4} & b_{5} \\ 0 & 0 & 0 & 3 a_{1} + a_{2} & b_{3} & b_{4} \\ 0 & 0 & 0 & 0 & 4 a_{1} + a_{2} & b_{3} \\ 0 & 0 & 0 & 0 & 0 & 5 a_{1} + a_{2}\end{pmatrix}}
{B=\begin{pmatrix}c_{1} & c_{2} \\ 0 & a_{1} + a_{2} + c_{1}\end{pmatrix}}
{U=\begin{pmatrix}u_{1} & u_{2} \\ v_{1} & v_{2} \\ 0 & - u_{1} + v_{1} \\ 0 & 0 \\ 0 & 0 \\ 0 & 0\end{pmatrix}}
{V=\begin{pmatrix}0 & 0 & 0 & 0 & 0 & w_{6} \\ 0 & 0 & 0 & 0 & 0 & 0\end{pmatrix}}
\derivcase{case:5-2-6}{F_{5,2}^{(6)}}
{A=\begin{pmatrix}a_{1} & a_{2} & a_{3} & a_{4} & a_{5} & a_{6} \\ 0 & a_{1} + a_{2} & b_{3} & b_{4} & b_{5} & b_{6} \\ 0 & 0 & 2 a_{1} + a_{2} & b_{3} & b_{4} & b_{5} \\ 0 & 0 & 0 & 3 a_{1} + a_{2} & b_{3} & b_{4} \\ 0 & 0 & 0 & 0 & 4 a_{1} + a_{2} & b_{3} \\ 0 & 0 & 0 & 0 & 0 & 5 a_{1} + a_{2}\end{pmatrix}}
{B=\begin{pmatrix}\tfrac{5 a_{1} + a_{2}}{2} & c_{2} \\ 0 & \tfrac{7 a_{1} + 3 a_{2}}{2}\end{pmatrix}}
{U=\begin{pmatrix}u_{1} & u_{2} \\ u_{1} + w_{5} & v_{2} \\ 0 & w_{5} \\ 0 & 0 \\ 0 & 0 \\ 0 & 0\end{pmatrix}}
{V=\begin{pmatrix}0 & 0 & 0 & 0 & w_{5} & w_{6} \\ 0 & 0 & 0 & 0 & 0 & u_{1} + w_{5}\end{pmatrix}}
\derivcase{case:5-2-7}{F_{5,2}^{(7)}}
{A=\begin{pmatrix}a_{1} & a_{2} & a_{3} & a_{4} & a_{5} & a_{6} \\ 0 & 3 a_{1} & b_{3} & b_{4} & b_{5} & b_{6} \\ 0 & 0 & 4 a_{1} & b_{3} & b_{4} & - a_{3} + b_{5} \\ 0 & 0 & 0 & 5 a_{1} & b_{3} & a_{2} + b_{4} \\ 0 & 0 & 0 & 0 & 6 a_{1} & b_{3} \\ 0 & 0 & 0 & 0 & 0 & 7 a_{1}\end{pmatrix}}
{B=\begin{pmatrix}c_{1} & c_{2} \\ 0 & a_{1} + c_{1}\end{pmatrix}}
{U=\begin{pmatrix}u_{1} & u_{2} \\ v_{1} & v_{2} \\ 0 & v_{1} \\ 0 & 0 \\ 0 & 0 \\ 0 & 0\end{pmatrix}}
{V=\begin{pmatrix}0 & 0 & 0 & 0 & w_{5} & w_{6} \\ 0 & 0 & 0 & 0 & 0 & w_{5}\end{pmatrix}}
\derivcase{case:5-2-8}{F_{5,2}^{(8)}}
{A=\begin{pmatrix}a_{1} & a_{2} & a_{3} & a_{4} & a_{5} & a_{6} \\ 0 & 3 a_{1} & c_{2} & - a_{2} & b_{5} & b_{6} \\ 0 & 0 & 4 a_{1} & c_{2} & - a_{2} & - a_{3} + b_{5} \\ 0 & 0 & 0 & 5 a_{1} & c_{2} & 0 \\ 0 & 0 & 0 & 0 & 6 a_{1} & c_{2} \\ 0 & 0 & 0 & 0 & 0 & 7 a_{1}\end{pmatrix}}
{B=\begin{pmatrix}\tfrac{5 a_{1}}{2} & c_{2} \\ 0 & \tfrac{7 a_{1}}{2}\end{pmatrix}}
{U=\begin{pmatrix}0 & u_{2} \\ 0 & 0 \\ 0 & 0 \\ 0 & 0 \\ 0 & 0 \\ 0 & 0\end{pmatrix}}
{V=\begin{pmatrix}0 & 0 & 0 & - u_{2} & w_{5} & w_{6} \\ 0 & 0 & 0 & 0 & - \tfrac{u_{2}}{2} & w_{5}\end{pmatrix}}
\derivcase{case:5-2-9}{F_{5,2}^{(9)}}
{A=\begin{pmatrix}a_{1} & a_{2} & a_{3} & a_{4} & a_{5} & a_{6} \\ 0 & 3 a_{1} & c_{2} & b_{4} & b_{5} & b_{6} \\ 0 & 0 & 4 a_{1} & c_{2} & b_{4} & - a_{3} + b_{5} \\ 0 & 0 & 0 & 5 a_{1} & c_{2} & a_{2} + b_{4} \\ 0 & 0 & 0 & 0 & 6 a_{1} & c_{2} \\ 0 & 0 & 0 & 0 & 0 & 7 a_{1}\end{pmatrix}}
{B=\begin{pmatrix}3 a_{1} & c_{2} \\ 0 & 4 a_{1}\end{pmatrix}}
{U=\begin{pmatrix}u_{1} & u_{2} \\ 0 & 0 \\ 0 & 0 \\ 0 & 0 \\ 0 & 0 \\ 0 & 0\end{pmatrix}}
{V=\begin{pmatrix}0 & 0 & 0 & - \tfrac{3 u_{1}}{2} & w_{5} & w_{6} \\ 0 & 0 & 0 & 0 & - \tfrac{u_{1}}{2} & \tfrac{u_{2} + 2 w_{5}}{2}\end{pmatrix}}
\derivcase{case:5-2-10}{F_{5,2}^{(10)}}
{A=\begin{pmatrix}a_{1} & a_{2} & a_{3} & a_{4} & a_{5} & a_{6} \\ 0 & 3 a_{1} & b_{3} & b_{4} & b_{5} & b_{6} \\ 0 & 0 & 4 a_{1} & b_{3} & b_{4} & - a_{3} + b_{5} \\ 0 & 0 & 0 & 5 a_{1} & b_{3} & a_{2} + b_{4} \\ 0 & 0 & 0 & 0 & 6 a_{1} & b_{3} \\ 0 & 0 & 0 & 0 & 0 & 7 a_{1}\end{pmatrix}}
{B=\begin{pmatrix}\tfrac{7 a_{1}}{2} & c_{2} \\ 0 & \tfrac{9 a_{1}}{2}\end{pmatrix}}
{U=\begin{pmatrix}u_{1} & u_{2} \\ 0 & v_{2} \\ 0 & 0 \\ 0 & 0 \\ 0 & 0 \\ 0 & 0\end{pmatrix}}
{V=\begin{pmatrix}0 & 0 & 0 & 0 & w_{5} & w_{6} \\ 0 & 0 & 0 & 0 & 0 & u_{1} + w_{5}\end{pmatrix}}
\derivcase{case:5-2-11}{F_{5,2}^{(11)}}
{A=\begin{pmatrix}a_{1} & 2 a_{1} & a_{3} & a_{4} & a_{5} & a_{6} \\ 0 & 3 a_{1} & b_{3} & b_{4} & b_{5} & b_{6} \\ 0 & 0 & 4 a_{1} & b_{3} & b_{4} & - a_{3} + b_{5} \\ 0 & 0 & 0 & 5 a_{1} & b_{3} & 2 a_{1} + b_{4} \\ 0 & 0 & 0 & 0 & 6 a_{1} & b_{3} \\ 0 & 0 & 0 & 0 & 0 & 7 a_{1}\end{pmatrix}}
{B=\begin{pmatrix}c_{1} & c_{2} \\ 0 & 3 a_{1} + c_{1}\end{pmatrix}}
{U=\begin{pmatrix}u_{1} & u_{2} \\ v_{1} & v_{2} \\ 0 & - u_{1} + v_{1} \\ 0 & 0 \\ 0 & 0 \\ 0 & 0\end{pmatrix}}
{V=\begin{pmatrix}0 & 0 & 0 & 0 & 0 & w_{6} \\ 0 & 0 & 0 & 0 & 0 & 0\end{pmatrix}}
\derivcase{case:5-2-12}{F_{5,2}^{(12)}}
{A=\begin{pmatrix}a_{1} & 2 a_{1} & a_{3} & a_{4} & a_{5} & a_{6} \\ 0 & 3 a_{1} & b_{3} & b_{4} & b_{5} & b_{6} \\ 0 & 0 & 4 a_{1} & b_{3} & b_{4} & - a_{3} + b_{5} \\ 0 & 0 & 0 & 5 a_{1} & b_{3} & 2 a_{1} + b_{4} \\ 0 & 0 & 0 & 0 & 6 a_{1} & b_{3} \\ 0 & 0 & 0 & 0 & 0 & 7 a_{1}\end{pmatrix}}
{B=\begin{pmatrix}\tfrac{7 a_{1}}{2} & c_{2} \\ 0 & \tfrac{13 a_{1}}{2}\end{pmatrix}}
{U=\begin{pmatrix}u_{1} & u_{2} \\ u_{1} + w_{5} & v_{2} \\ 0 & w_{5} \\ 0 & 0 \\ 0 & 0 \\ 0 & 0\end{pmatrix}}
{V=\begin{pmatrix}0 & 0 & 0 & 0 & w_{5} & w_{6} \\ 0 & 0 & 0 & 0 & 0 & u_{1} + w_{5}\end{pmatrix}}
\derivcase{case:5-2-13}{F_{5,2}^{(13)}}
{A=\begin{pmatrix}a_{1} & - a_{1} & a_{3} & a_{4} & a_{5} & a_{6} \\ 0 & 2 a_{1} & b_{3} & a_{1} & b_{5} & b_{6} \\ 0 & 0 & 3 a_{1} & b_{3} & a_{1} - a_{3} & - a_{4} + a_{5} + b_{5} \\ 0 & 0 & 0 & 4 a_{1} & - a_{1} + b_{3} & a_{1} - a_{3} - a_{4} \\ 0 & 0 & 0 & 0 & 5 a_{1} & - 2 a_{1} + a_{3} + b_{3} \\ 0 & 0 & 0 & 0 & 0 & 7 a_{1}\end{pmatrix}}
{B=\begin{pmatrix}c_{1} & c_{2} \\ 0 & a_{1} + c_{1}\end{pmatrix}}
{U=\begin{pmatrix}u_{1} & u_{2} \\ v_{1} & v_{2} \\ 0 & v_{1} \\ 0 & 0 \\ 0 & 0 \\ 0 & 0\end{pmatrix}}
{V=\begin{pmatrix}0 & 0 & 0 & 0 & 0 & w_{6} \\ 0 & 0 & 0 & 0 & 0 & 0\end{pmatrix}}
\derivcase{case:5-2-14}{F_{5,2}^{(14)}}
{A=\begin{pmatrix}a_{1} & - a_{1} & a_{3} & a_{4} & a_{5} & a_{6} \\ 0 & 2 a_{1} & b_{3} & a_{1} & b_{5} & b_{6} \\ 0 & 0 & 3 a_{1} & b_{3} & a_{1} - a_{3} & - a_{4} + a_{5} + b_{5} \\ 0 & 0 & 0 & 4 a_{1} & - a_{1} + b_{3} & a_{1} - a_{3} - a_{4} \\ 0 & 0 & 0 & 0 & 5 a_{1} & - 2 a_{1} + a_{3} + b_{3} \\ 0 & 0 & 0 & 0 & 0 & 7 a_{1}\end{pmatrix}}
{B=\begin{pmatrix}\tfrac{7 a_{1}}{2} & c_{2} \\ 0 & \tfrac{9 a_{1}}{2}\end{pmatrix}}
{U=\begin{pmatrix}u_{1} & u_{2} \\ w_{5} & v_{2} \\ 0 & w_{5} \\ 0 & 0 \\ 0 & 0 \\ 0 & 0\end{pmatrix}}
{V=\begin{pmatrix}0 & 0 & 0 & 0 & w_{5} & w_{6} \\ 0 & 0 & 0 & 0 & 0 & u_{1} + w_{5}\end{pmatrix}}
\derivcase{case:5-2-15}{F_{5,2}^{(15)}}
{A=\begin{pmatrix}0 & 0 & a_{3} & a_{4} & a_{5} & a_{6} \\ 0 & 0 & b_{3} & 0 & b_{5} & b_{6} \\ 0 & 0 & 0 & b_{3} & - a_{3} & - a_{4} + a_{5} + b_{5} \\ 0 & 0 & 0 & 0 & b_{3} & - a_{3} - a_{4} \\ 0 & 0 & 0 & 0 & 0 & a_{3} + b_{3} \\ 0 & 0 & 0 & 0 & 0 & 0\end{pmatrix}}
{B=\begin{pmatrix}c_{1} & c_{2} \\ 0 & c_{1}\end{pmatrix}}
{U=\begin{pmatrix}u_{1} & u_{2} \\ v_{1} & v_{2} \\ 0 & - u_{1} + v_{1} \\ 0 & 0 \\ 0 & 0 \\ 0 & 0\end{pmatrix}}
{V=\begin{pmatrix}0 & 0 & 0 & 0 & 0 & w_{6} \\ 0 & 0 & 0 & 0 & 0 & 0\end{pmatrix}}
\derivcase{case:5-2-16}{F_{5,2}^{(16)}}
{A=\begin{pmatrix}0 & 0 & a_{3} & a_{4} & a_{5} & a_{6} \\ 0 & 0 & b_{3} & 0 & b_{5} & b_{6} \\ 0 & 0 & 0 & b_{3} & - a_{3} & - a_{4} + a_{5} + b_{5} \\ 0 & 0 & 0 & 0 & b_{3} & - a_{3} - a_{4} \\ 0 & 0 & 0 & 0 & 0 & a_{3} + b_{3} \\ 0 & 0 & 0 & 0 & 0 & 0\end{pmatrix}}
{B=\begin{pmatrix}0 & c_{2} \\ 0 & 0\end{pmatrix}}
{U=\begin{pmatrix}u_{1} & u_{2} \\ u_{1} + 2 w_{5} & v_{2} \\ 0 & 2 w_{5} \\ 0 & 0 \\ 0 & 0 \\ 0 & 0\end{pmatrix}}
{V=\begin{pmatrix}0 & 0 & 0 & 0 & w_{5} & w_{6} \\ 0 & 0 & 0 & 0 & 0 & u_{1} + w_{5}\end{pmatrix}}
\derivcase{case:5-2-17}{F_{5,2}^{(17)}}
{A=\begin{pmatrix}a_{1} & a_{2} & a_{3} & a_{4} & a_{5} & a_{6} \\ 0 & a_{1} - a_{2} & b_{3} & 0 & b_{5} & b_{6} \\ 0 & 0 & 2 a_{1} - a_{2} & b_{3} & 0 & a_{5} + b_{5} \\ 0 & 0 & 0 & 3 a_{1} - a_{2} & b_{3} & - a_{4} \\ 0 & 0 & 0 & 0 & 4 a_{1} - a_{2} & a_{3} + b_{3} \\ 0 & 0 & 0 & 0 & 0 & 5 a_{1} - 2 a_{2}\end{pmatrix}}
{B=\begin{pmatrix}c_{1} & c_{2} \\ 0 & a_{1} + c_{1}\end{pmatrix}}
{U=\begin{pmatrix}u_{1} & u_{2} \\ v_{1} & v_{2} \\ 0 & v_{1} \\ 0 & 0 \\ 0 & 0 \\ 0 & 0\end{pmatrix}}
{V=\begin{pmatrix}0 & 0 & 0 & 0 & 0 & w_{6} \\ 0 & 0 & 0 & 0 & 0 & 0\end{pmatrix}}
\derivcase{case:5-2-18}{F_{5,2}^{(18)}}
{A=\begin{pmatrix}a_{1} & a_{2} & a_{3} & a_{4} & a_{5} & a_{6} \\ 0 & a_{1} - a_{2} & b_{3} & 0 & b_{5} & b_{6} \\ 0 & 0 & 2 a_{1} - a_{2} & b_{3} & 0 & a_{5} + b_{5} \\ 0 & 0 & 0 & 3 a_{1} - a_{2} & b_{3} & - a_{4} \\ 0 & 0 & 0 & 0 & 4 a_{1} - a_{2} & a_{3} + b_{3} \\ 0 & 0 & 0 & 0 & 0 & 5 a_{1} - 2 a_{2}\end{pmatrix}}
{B=\begin{pmatrix}\tfrac{5 a_{1} - 2 a_{2}}{2} & c_{2} \\ 0 & \tfrac{7 a_{1} - 2 a_{2}}{2}\end{pmatrix}}
{U=\begin{pmatrix}u_{1} & u_{2} \\ w_{5} & v_{2} \\ 0 & w_{5} \\ 0 & 0 \\ 0 & 0 \\ 0 & 0\end{pmatrix}}
{V=\begin{pmatrix}0 & 0 & 0 & 0 & w_{5} & w_{6} \\ 0 & 0 & 0 & 0 & 0 & u_{1} + w_{5}\end{pmatrix}}
\derivcase{case:5-2-19}{F_{5,2}^{(19)}}
{A=\begin{pmatrix}a_{1} & 0 & a_{3} & a_{4} & a_{5} & a_{6} \\ 0 & a_{1} & b_{3} & 0 & b_{5} & b_{6} \\ 0 & 0 & 2 a_{1} & b_{3} & 0 & a_{5} + b_{5} \\ 0 & 0 & 0 & 3 a_{1} & b_{3} & - a_{4} \\ 0 & 0 & 0 & 0 & 4 a_{1} & a_{3} + b_{3} \\ 0 & 0 & 0 & 0 & 0 & 5 a_{1}\end{pmatrix}}
{B=\begin{pmatrix}c_{1} & c_{2} \\ 0 & a_{1} + c_{1}\end{pmatrix}}
{U=\begin{pmatrix}u_{1} & u_{2} \\ v_{1} & v_{2} \\ 0 & - u_{1} + v_{1} \\ 0 & 0 \\ 0 & 0 \\ 0 & 0\end{pmatrix}}
{V=\begin{pmatrix}0 & 0 & 0 & 0 & 0 & w_{6} \\ 0 & 0 & 0 & 0 & 0 & 0\end{pmatrix}}
\derivcase{case:5-2-20}{F_{5,2}^{(20)}}
{A=\begin{pmatrix}a_{1} & 0 & a_{3} & a_{4} & a_{5} & a_{6} \\ 0 & a_{1} & b_{3} & 0 & b_{5} & b_{6} \\ 0 & 0 & 2 a_{1} & b_{3} & 0 & a_{5} + b_{5} \\ 0 & 0 & 0 & 3 a_{1} & b_{3} & - a_{4} \\ 0 & 0 & 0 & 0 & 4 a_{1} & a_{3} + b_{3} \\ 0 & 0 & 0 & 0 & 0 & 5 a_{1}\end{pmatrix}}
{B=\begin{pmatrix}\tfrac{5 a_{1}}{2} & c_{2} \\ 0 & \tfrac{7 a_{1}}{2}\end{pmatrix}}
{U=\begin{pmatrix}u_{1} & u_{2} \\ u_{1} + 2 w_{5} & v_{2} \\ 0 & 2 w_{5} \\ 0 & 0 \\ 0 & 0 \\ 0 & 0\end{pmatrix}}
{V=\begin{pmatrix}0 & 0 & 0 & 0 & w_{5} & w_{6} \\ 0 & 0 & 0 & 0 & 0 & u_{1} + w_{5}\end{pmatrix}}
\derivcase{case:5-2-21}{F_{5,2}^{(21)}}
{A=\begin{pmatrix}a_{1} & a_{2} & a_{3} & a_{4} & a_{5} & a_{6} \\ 0 & a_{1} - a_{2} & b_{3} & 0 & b_{5} & b_{6} \\ 0 & 0 & 2 a_{1} - a_{2} & b_{3} & 0 & a_{5} + b_{5} \\ 0 & 0 & 0 & 3 a_{1} - a_{2} & b_{3} & - a_{4} \\ 0 & 0 & 0 & 0 & 4 a_{1} - a_{2} & a_{3} + b_{3} \\ 0 & 0 & 0 & 0 & 0 & 5 a_{1} - 2 a_{2}\end{pmatrix}}
{B=\begin{pmatrix}c_{1} & c_{2} \\ 0 & a_{1} - a_{2} + c_{1}\end{pmatrix}}
{U=\begin{pmatrix}u_{1} & u_{2} \\ v_{1} & v_{2} \\ 0 & u_{1} + v_{1} \\ 0 & 0 \\ 0 & 0 \\ 0 & 0\end{pmatrix}}
{V=\begin{pmatrix}0 & 0 & 0 & 0 & w_{5} & w_{6} \\ 0 & 0 & 0 & 0 & 0 & w_{5}\end{pmatrix}}
\derivcase{case:5-2-22}{F_{5,2}^{(22)}}
{A=\begin{pmatrix}a_{1} & a_{2} & a_{3} & a_{4} & a_{5} & a_{6} \\ 0 & a_{1} - a_{2} & - a_{3} + c_{2} & 0 & b_{5} & b_{6} \\ 0 & 0 & 2 a_{1} - a_{2} & - a_{3} + c_{2} & 0 & a_{5} + b_{5} \\ 0 & 0 & 0 & 3 a_{1} - a_{2} & - a_{3} + c_{2} & - a_{4} \\ 0 & 0 & 0 & 0 & 4 a_{1} - a_{2} & c_{2} \\ 0 & 0 & 0 & 0 & 0 & 5 a_{1} - 2 a_{2}\end{pmatrix}}
{B=\begin{pmatrix}\tfrac{4 a_{1} - a_{2}}{2} & c_{2} \\ 0 & \tfrac{3 \left(2 a_{1} - a_{2}\right)}{2}\end{pmatrix}}
{U=\begin{pmatrix}0 & u_{2} \\ 0 & - u_{2} \\ 0 & 0 \\ 0 & 0 \\ 0 & 0 \\ 0 & 0\end{pmatrix}}
{V=\begin{pmatrix}0 & 0 & 0 & 0 & w_{5} & w_{6} \\ 0 & 0 & 0 & 0 & 0 & \tfrac{u_{2} + 2 w_{5}}{2}\end{pmatrix}}
\derivcase{case:5-2-23}{F_{5,2}^{(23)}}
{A=\begin{pmatrix}a_{1} & a_{2} & a_{3} & a_{4} & a_{5} & a_{6} \\ 0 & a_{1} - a_{2} & b_{3} & 0 & b_{5} & b_{6} \\ 0 & 0 & 2 a_{1} - a_{2} & b_{3} & 0 & a_{5} + b_{5} \\ 0 & 0 & 0 & 3 a_{1} - a_{2} & b_{3} & - a_{4} \\ 0 & 0 & 0 & 0 & 4 a_{1} - a_{2} & a_{3} + b_{3} \\ 0 & 0 & 0 & 0 & 0 & 5 a_{1} - 2 a_{2}\end{pmatrix}}
{B=\begin{pmatrix}\tfrac{5 a_{1} - 2 a_{2}}{2} & c_{2} \\ 0 & \tfrac{7 a_{1} - 4 a_{2}}{2}\end{pmatrix}}
{U=\begin{pmatrix}u_{1} & u_{2} \\ - u_{1} & v_{2} \\ 0 & 0 \\ 0 & 0 \\ 0 & 0 \\ 0 & 0\end{pmatrix}}
{V=\begin{pmatrix}0 & 0 & 0 & 0 & w_{5} & w_{6} \\ 0 & 0 & 0 & 0 & 0 & u_{1} + w_{5}\end{pmatrix}}
\derivcase{case:5-2-24}{F_{5,2}^{(24)}}
{A=\begin{pmatrix}a_{1} & a_{2} & a_{3} & a_{4} & a_{5} & a_{6} \\ 0 & a_{1} - a_{2} & b_{3} & - a_{1} & b_{5} & b_{6} \\ 0 & 0 & 2 a_{1} - a_{2} & b_{3} & - a_{1} & - a_{3} + a_{5} + b_{5} \\ 0 & 0 & 0 & 3 a_{1} - a_{2} & b_{3} & - a_{1} + a_{2} - a_{4} \\ 0 & 0 & 0 & 0 & 4 a_{1} - a_{2} & a_{3} + b_{3} \\ 0 & 0 & 0 & 0 & 0 & 5 a_{1} - 2 a_{2}\end{pmatrix}}
{B=\begin{pmatrix}c_{1} & c_{2} \\ 0 & a_{1} + c_{1}\end{pmatrix}}
{U=\begin{pmatrix}u_{1} & u_{2} \\ v_{1} & v_{2} \\ 0 & v_{1} \\ 0 & 0 \\ 0 & 0 \\ 0 & 0\end{pmatrix}}
{V=\begin{pmatrix}0 & 0 & 0 & 0 & 0 & w_{6} \\ 0 & 0 & 0 & 0 & 0 & 0\end{pmatrix}}
\derivcase{case:5-2-25}{F_{5,2}^{(25)}}
{A=\begin{pmatrix}a_{1} & a_{2} & a_{3} & a_{4} & a_{5} & a_{6} \\ 0 & a_{1} - a_{2} & b_{3} & - a_{1} & b_{5} & b_{6} \\ 0 & 0 & 2 a_{1} - a_{2} & b_{3} & - a_{1} & - a_{3} + a_{5} + b_{5} \\ 0 & 0 & 0 & 3 a_{1} - a_{2} & b_{3} & - a_{1} + a_{2} - a_{4} \\ 0 & 0 & 0 & 0 & 4 a_{1} - a_{2} & a_{3} + b_{3} \\ 0 & 0 & 0 & 0 & 0 & 5 a_{1} - 2 a_{2}\end{pmatrix}}
{B=\begin{pmatrix}\tfrac{5 a_{1} - 2 a_{2}}{2} & c_{2} \\ 0 & \tfrac{7 a_{1} - 2 a_{2}}{2}\end{pmatrix}}
{U=\begin{pmatrix}u_{1} & u_{2} \\ w_{5} & v_{2} \\ 0 & w_{5} \\ 0 & 0 \\ 0 & 0 \\ 0 & 0\end{pmatrix}}
{V=\begin{pmatrix}0 & 0 & 0 & 0 & w_{5} & w_{6} \\ 0 & 0 & 0 & 0 & 0 & u_{1} + w_{5}\end{pmatrix}}
\derivcase{case:5-2-26}{F_{5,2}^{(26)}}
{A=\begin{pmatrix}a_{1} & 0 & a_{3} & a_{4} & a_{5} & a_{6} \\ 0 & a_{1} & b_{3} & - a_{1} & b_{5} & b_{6} \\ 0 & 0 & 2 a_{1} & b_{3} & - a_{1} & - a_{3} + a_{5} + b_{5} \\ 0 & 0 & 0 & 3 a_{1} & b_{3} & - a_{1} - a_{4} \\ 0 & 0 & 0 & 0 & 4 a_{1} & a_{3} + b_{3} \\ 0 & 0 & 0 & 0 & 0 & 5 a_{1}\end{pmatrix}}
{B=\begin{pmatrix}c_{1} & c_{2} \\ 0 & a_{1} + c_{1}\end{pmatrix}}
{U=\begin{pmatrix}u_{1} & u_{2} \\ v_{1} & v_{2} \\ 0 & - u_{1} + v_{1} \\ 0 & 0 \\ 0 & 0 \\ 0 & 0\end{pmatrix}}
{V=\begin{pmatrix}0 & 0 & 0 & 0 & 0 & w_{6} \\ 0 & 0 & 0 & 0 & 0 & 0\end{pmatrix}}
\derivcase{case:5-2-27}{F_{5,2}^{(27)}}
{A=\begin{pmatrix}a_{1} & 0 & a_{3} & a_{4} & a_{5} & a_{6} \\ 0 & a_{1} & b_{3} & - a_{1} & b_{5} & b_{6} \\ 0 & 0 & 2 a_{1} & b_{3} & - a_{1} & - a_{3} + a_{5} + b_{5} \\ 0 & 0 & 0 & 3 a_{1} & b_{3} & - a_{1} - a_{4} \\ 0 & 0 & 0 & 0 & 4 a_{1} & a_{3} + b_{3} \\ 0 & 0 & 0 & 0 & 0 & 5 a_{1}\end{pmatrix}}
{B=\begin{pmatrix}\tfrac{5 a_{1}}{2} & c_{2} \\ 0 & \tfrac{7 a_{1}}{2}\end{pmatrix}}
{U=\begin{pmatrix}u_{1} & u_{2} \\ u_{1} + 2 w_{5} & v_{2} \\ 0 & 2 w_{5} \\ 0 & 0 \\ 0 & 0 \\ 0 & 0\end{pmatrix}}
{V=\begin{pmatrix}0 & 0 & 0 & 0 & w_{5} & w_{6} \\ 0 & 0 & 0 & 0 & 0 & u_{1} + w_{5}\end{pmatrix}}
\derivcase{case:5-2-28}{F_{5,2}^{(28)}}
{A=\begin{pmatrix}a_{1} & a_{2} & a_{3} & a_{4} & a_{5} & a_{6} \\ 0 & a_{1} - a_{2} & b_{3} & - a_{1} & b_{5} & b_{6} \\ 0 & 0 & 2 a_{1} - a_{2} & b_{3} & - a_{1} & - a_{3} + a_{5} + b_{5} \\ 0 & 0 & 0 & 3 a_{1} - a_{2} & b_{3} & - a_{1} + a_{2} - a_{4} \\ 0 & 0 & 0 & 0 & 4 a_{1} - a_{2} & a_{3} + b_{3} \\ 0 & 0 & 0 & 0 & 0 & 5 a_{1} - 2 a_{2}\end{pmatrix}}
{B=\begin{pmatrix}c_{1} & c_{2} \\ 0 & a_{1} - a_{2} + c_{1}\end{pmatrix}}
{U=\begin{pmatrix}u_{1} & u_{2} \\ v_{1} & v_{2} \\ 0 & u_{1} + v_{1} \\ 0 & 0 \\ 0 & 0 \\ 0 & 0\end{pmatrix}}
{V=\begin{pmatrix}0 & 0 & 0 & 0 & w_{5} & w_{6} \\ 0 & 0 & 0 & 0 & 0 & w_{5}\end{pmatrix}}
\derivcase{case:5-2-29}{F_{5,2}^{(29)}}
{A=\begin{pmatrix}a_{1} & a_{2} & a_{3} & a_{4} & a_{5} & a_{6} \\ 0 & a_{1} - a_{2} & - a_{3} + c_{2} & - a_{1} & b_{5} & b_{6} \\ 0 & 0 & 2 a_{1} - a_{2} & - a_{3} + c_{2} & - a_{1} & - a_{3} + a_{5} + b_{5} \\ 0 & 0 & 0 & 3 a_{1} - a_{2} & - a_{3} + c_{2} & - a_{1} + a_{2} - a_{4} \\ 0 & 0 & 0 & 0 & 4 a_{1} - a_{2} & c_{2} \\ 0 & 0 & 0 & 0 & 0 & 5 a_{1} - 2 a_{2}\end{pmatrix}}
{B=\begin{pmatrix}\tfrac{4 a_{1} - a_{2}}{2} & c_{2} \\ 0 & \tfrac{3 \left(2 a_{1} - a_{2}\right)}{2}\end{pmatrix}}
{U=\begin{pmatrix}0 & u_{2} \\ 0 & - u_{2} \\ 0 & 0 \\ 0 & 0 \\ 0 & 0 \\ 0 & 0\end{pmatrix}}
{V=\begin{pmatrix}0 & 0 & 0 & 0 & w_{5} & w_{6} \\ 0 & 0 & 0 & 0 & 0 & \tfrac{u_{2} + 2 w_{5}}{2}\end{pmatrix}}
\derivcase{case:5-2-30}{F_{5,2}^{(30)}}
{A=\begin{pmatrix}a_{1} & a_{2} & a_{3} & a_{4} & a_{5} & a_{6} \\ 0 & a_{1} - a_{2} & b_{3} & - a_{1} & b_{5} & b_{6} \\ 0 & 0 & 2 a_{1} - a_{2} & b_{3} & - a_{1} & - a_{3} + a_{5} + b_{5} \\ 0 & 0 & 0 & 3 a_{1} - a_{2} & b_{3} & - a_{1} + a_{2} - a_{4} \\ 0 & 0 & 0 & 0 & 4 a_{1} - a_{2} & a_{3} + b_{3} \\ 0 & 0 & 0 & 0 & 0 & 5 a_{1} - 2 a_{2}\end{pmatrix}}
{B=\begin{pmatrix}\tfrac{5 a_{1} - 2 a_{2}}{2} & c_{2} \\ 0 & \tfrac{7 a_{1} - 4 a_{2}}{2}\end{pmatrix}}
{U=\begin{pmatrix}u_{1} & u_{2} \\ - u_{1} & v_{2} \\ 0 & 0 \\ 0 & 0 \\ 0 & 0 \\ 0 & 0\end{pmatrix}}
{V=\begin{pmatrix}0 & 0 & 0 & 0 & w_{5} & w_{6} \\ 0 & 0 & 0 & 0 & 0 & u_{1} + w_{5}\end{pmatrix}}
\derivcase{case:5-2-31}{F_{5,2}^{(31)}}
{A=\begin{pmatrix}a_{1} & - \tfrac{a_{1}}{2} & a_{3} & a_{4} & a_{5} & a_{6} \\ 0 & 2 a_{1} & b_{3} & b_{4} & b_{5} & b_{6} \\ 0 & 0 & 3 a_{1} & b_{3} & - a_{3} + b_{4} & - a_{3} - a_{4} + b_{5} \\ 0 & 0 & 0 & 4 a_{1} & - \tfrac{a_{1} - 2 b_{3}}{2} & - \tfrac{a_{1} + 2 a_{3} - 2 b_{4}}{2} \\ 0 & 0 & 0 & 0 & 5 a_{1} & - a_{1} + b_{3} \\ 0 & 0 & 0 & 0 & 0 & 6 a_{1}\end{pmatrix}}
{B=\begin{pmatrix}c_{1} & c_{2} \\ 0 & 2 a_{1} + c_{1}\end{pmatrix}}
{U=\begin{pmatrix}u_{1} & u_{2} \\ v_{1} & v_{2} \\ 0 & 2 u_{1} + v_{1} \\ 0 & 0 \\ 0 & 0 \\ 0 & 0\end{pmatrix}}
{V=\begin{pmatrix}0 & 0 & 0 & 0 & 0 & w_{6} \\ 0 & 0 & 0 & 0 & 0 & 0\end{pmatrix}}
\derivcase{case:5-2-32}{F_{5,2}^{(32)}}
{A=\begin{pmatrix}a_{1} & - \tfrac{a_{1}}{2} & a_{3} & a_{4} & a_{5} & a_{6} \\ 0 & 2 a_{1} & b_{3} & b_{4} & b_{5} & b_{6} \\ 0 & 0 & 3 a_{1} & b_{3} & - a_{3} + b_{4} & - a_{3} - a_{4} + b_{5} \\ 0 & 0 & 0 & 4 a_{1} & - \tfrac{a_{1} - 2 b_{3}}{2} & - \tfrac{a_{1} + 2 a_{3} - 2 b_{4}}{2} \\ 0 & 0 & 0 & 0 & 5 a_{1} & - a_{1} + b_{3} \\ 0 & 0 & 0 & 0 & 0 & 6 a_{1}\end{pmatrix}}
{B=\begin{pmatrix}3 a_{1} & c_{2} \\ 0 & 5 a_{1}\end{pmatrix}}
{U=\begin{pmatrix}u_{1} & u_{2} \\ - 2 \left(u_{1} + w_{5}\right) & v_{2} \\ 0 & - 2 w_{5} \\ 0 & 0 \\ 0 & 0 \\ 0 & 0\end{pmatrix}}
{V=\begin{pmatrix}0 & 0 & 0 & 0 & w_{5} & w_{6} \\ 0 & 0 & 0 & 0 & 0 & u_{1} + w_{5}\end{pmatrix}}
\derivcase{case:5-2-33}{F_{5,2}^{(33)}}
{A=\begin{pmatrix}a_{1} & - a_{1} & a_{3} & a_{4} & a_{5} & a_{6} \\ 0 & 2 a_{1} & b_{3} & a_{1} & b_{5} & b_{6} \\ 0 & 0 & 3 a_{1} & b_{3} & a_{1} - a_{3} & - a_{4} + a_{5} + b_{5} \\ 0 & 0 & 0 & 4 a_{1} & - a_{1} + b_{3} & a_{1} - a_{3} - a_{4} \\ 0 & 0 & 0 & 0 & 5 a_{1} & - 2 a_{1} + a_{3} + b_{3} \\ 0 & 0 & 0 & 0 & 0 & 7 a_{1}\end{pmatrix}}
{B=\begin{pmatrix}c_{1} & c_{2} \\ 0 & 2 a_{1} + c_{1}\end{pmatrix}}
{U=\begin{pmatrix}u_{1} & u_{2} \\ v_{1} & v_{2} \\ 0 & u_{1} + v_{1} \\ 0 & 0 \\ 0 & 0 \\ 0 & 0\end{pmatrix}}
{V=\begin{pmatrix}0 & 0 & 0 & 0 & w_{5} & w_{6} \\ 0 & 0 & 0 & 0 & 0 & w_{5}\end{pmatrix}}
\derivcase{case:5-2-34}{F_{5,2}^{(34)}}
{A=\begin{pmatrix}a_{1} & - a_{1} & a_{3} & a_{4} & a_{5} & a_{6} \\ 0 & 2 a_{1} & 2 a_{1} - a_{3} + c_{2} & a_{1} & b_{5} & b_{6} \\ 0 & 0 & 3 a_{1} & 2 a_{1} - a_{3} + c_{2} & a_{1} - a_{3} & - a_{4} + a_{5} + b_{5} \\ 0 & 0 & 0 & 4 a_{1} & a_{1} - a_{3} + c_{2} & a_{1} - a_{3} - a_{4} \\ 0 & 0 & 0 & 0 & 5 a_{1} & c_{2} \\ 0 & 0 & 0 & 0 & 0 & 7 a_{1}\end{pmatrix}}
{B=\begin{pmatrix}\tfrac{5 a_{1}}{2} & c_{2} \\ 0 & \tfrac{9 a_{1}}{2}\end{pmatrix}}
{U=\begin{pmatrix}0 & u_{2} \\ 0 & - u_{2} \\ 0 & 0 \\ 0 & 0 \\ 0 & 0 \\ 0 & 0\end{pmatrix}}
{V=\begin{pmatrix}0 & 0 & 0 & 0 & w_{5} & w_{6} \\ 0 & 0 & 0 & 0 & 0 & \tfrac{u_{2} + 2 w_{5}}{2}\end{pmatrix}}
\derivcase{case:5-2-35}{F_{5,2}^{(35)}}
{A=\begin{pmatrix}a_{1} & - a_{1} & a_{3} & a_{4} & a_{5} & a_{6} \\ 0 & 2 a_{1} & b_{3} & a_{1} & b_{5} & b_{6} \\ 0 & 0 & 3 a_{1} & b_{3} & a_{1} - a_{3} & - a_{4} + a_{5} + b_{5} \\ 0 & 0 & 0 & 4 a_{1} & - a_{1} + b_{3} & a_{1} - a_{3} - a_{4} \\ 0 & 0 & 0 & 0 & 5 a_{1} & - 2 a_{1} + a_{3} + b_{3} \\ 0 & 0 & 0 & 0 & 0 & 7 a_{1}\end{pmatrix}}
{B=\begin{pmatrix}\tfrac{7 a_{1}}{2} & c_{2} \\ 0 & \tfrac{11 a_{1}}{2}\end{pmatrix}}
{U=\begin{pmatrix}u_{1} & u_{2} \\ - u_{1} & v_{2} \\ 0 & 0 \\ 0 & 0 \\ 0 & 0 \\ 0 & 0\end{pmatrix}}
{V=\begin{pmatrix}0 & 0 & 0 & 0 & w_{5} & w_{6} \\ 0 & 0 & 0 & 0 & 0 & u_{1} + w_{5}\end{pmatrix}}
\end{document}